\pdfoutput=1
\documentclass[11pt,intlimits]{article}

\usepackage{fancyhdr}
\usepackage{amssymb,amsmath,amsthm}
\usepackage{graphicx,float,wrapfig}
\usepackage{multirow}
\begin{document}
\title{Micronuclear Sequences Associated with Assembly Graphs}
 \author{Christeen Bisnath}
\maketitle
\begin{abstract}
\indent This paper investigates given an assembly graph, find the possible micronuclear sequences in terms of MDSs and IESs; which represent Hamiltonian polygonal paths. We will consider the orientations of the assembly graph and Hamiltonian polygonal path. To obtain a micronuclear sequence of this path, we will compare orientation of the Hamiltonian polygonal path with respect to the orientation of the assembly graph. We concentrate on two assembly graphs with two 4-valent vertices. These Hamiltonian polygonal paths influence the type of smoothings performed. There are two types of smoothings parallel and non-parallel smoothing. We consider a micronuclear sequence after one smoothing of a vertex. 
\end{abstract}
\section{Introduction}
%Detail intro., give references, including the association of hamiltonian polygonal path purpose, define transversal, fix second last sentence, define transversal\\
% A Hamiltonian path visits each vertex exactly once, polygonal path visits each vertex in a perpendicular direction. The micronuclear sequences include MDS (macronuclear destined segment) and IES (internal eliminated segment). \\
\indent During the macronuclear development of ciliates \textit{(Stylonychia or Oxytricha)}, DNA recombination occurs. Each nucleotide sequence  of the macronucleus may be in non-consecutive order in the micronucleus as macronuclear destined sequence (MDSs) separated by non-coding DNA. During macronuclear development the non-coding fragments (internal eliminated sequence, IESs) that interrupt MDSs in the micronucleus are excised. However, the order of MDSs in the micronucleus may not be consecutive, and may be reversed relative to the others \cite{1}. The rearrangement of a micronuclear sequence can be seen as MDSs and IESs. The assembly graphs represent these micronuclear genes. An assembly graph is a finite connected graph where all vertices are rigid vertices of valency 1 or 4 \cite{1}. We will focus on assembly graphs  that have Hamiltonian polygonal paths that represent a single gene combination in the micronucleus, such graphs are called \textit{realizable graphs}. A polygonal path is a path in which every two consecutive edges are neighbors with respect to their common vertex or in other words, makes a "90 degree" turn at every rigid vertex \cite{1}. A Hamiltonian polygonal path is a path in which each 4-valent is visited exactly once, makes a "90 degree turn" at each rigid vertex and no edges are repeated. A transverse path is a path in which every two consecutive edges are not neighbors. The given assembly graph (transversal) corresponds  to the micronuclear DNA segment before recombination. A transverse graph (transversal) is where a single path contains every edge is visited and none is repeated, which are called simple assembly graphs \cite{3}. We are interested in what types of rearrangements and how many distinct genes can be encoded by an assembly graph \cite{4}.\\
%explain these are assembly graphs that have a single hamiltonian polygonal path (realizable graphs), given the hamiltonian polygonal paths figure out how to write the mic sequence using the labellings of the assembly graph
%explain these are assembly graphs that have a single hamiltonian polygonal path (realizable graphs), given the hamiltonian polygonal paths figure out how to write the mic sequence using the labellings of the assembly graph
\section{Forming Micronuclear Sequences}
\subsection{Notation for Hamiltonian Polygonal Paths and Micronuclear Sequences}
\indent Given a transversal with possible Hamiltonian polygonal paths (HPPs), we will write possible DNA recombinations of micronuclear sequences. In order to understand Hamiltonian polygonal paths and micronuclear sequences, the following will be some notation to help describe them.
%Now, given a transversal graph, to find the possible Hamiltonian polygonal paths, pick an end of the assembly graph as the reference point. This will stay as the original path of the original assembly graph (transversal). % To obtain the micronulcear sequence, given a Hamiltonian polygonal path, use the definitions stated previously; each path should be perpendicular at all vertices and non repeating edges and vertices. Use the minimal amount of edges at a time when finding the paths. The sequence is read following the order of the original path.To write the micronuclear sequence, start comparing the Hamiltonian polygonal path to the original path. The sequence is read following the order of the original path In order to write the  micronuclear sequence, let's introduce some notations, to better describe the assembly graphs and micronuclear sequences.\\
\begin{itemize}
\item $n$ represent 4-valent vertices\\
\item $e_0, \ldots,e_{2n}$ represent edges, are enumerated by the way they are encountered by the transverse path, this fixes an orientation of the edge\\
\item $v_0,\dots,v_{n+1}$ represent vertices of assembly graph\\
\item $\Gamma={\gamma_1,\dots,\gamma_{n+1}}$ where $\Gamma$ is one HPP and  $\gamma_1, \gamma_{n+1}$ are ending edges and $\gamma_2,\dots, \gamma_n$ are full edges of G and $\Gamma$, let G assembly graph be $n$ 4-valent vertices
\item $P$= $\gamma_2,\dots, \gamma_n$\\
\item $I_0,\ldots,I_{n+1}$ represent IES sequences, for n is number of vertices in assembly graph\\
\item $M_1,\ldots,M_{n+1}$ represent MDS sequences, for n is number of vertices in assembly graph\\
\item $\overline{M_1},\ldots, \overline{M_{n+1}}$ represent inverted MDS sequences, for n is number of vertices in assembly graph\\
\item - - - - - - - - - dashed lines represent the Hamiltonian polygonal path\\
\item \line(1,0){67} \hspace{1mm} solid line represent the transversal path\\
\item \noindent A= incoming edges, B= outgoing edges; of transversal path orientation\\
We indicate with A an incoming edge incident to a vertex and with B outgoing edge.
\end{itemize}

\indent Now, with a given assembly graph and Hamiltonian polygonal path without ending edges (P), consider the first and last ending edges of the path, ($\gamma_1, \gamma_{n+1}$) and vertices of the path. For P= $\gamma_2$,$\dots$,$\gamma_{n}$, which are the complete edges from the assembly graph. Also, HPP= $w_0$ $\gamma_1$,$\dots$, $w_{n-1}$$\gamma_{n+1}$; the order of encountering these edges defines an orientation of the HPP. Label the incoming and outgoing edges that are neighbors of $\gamma_1$ incident to $v_{i_0}$ with 1A and 1B respectively. Then, label the incoming and outgoing edges that are neighbors of $\gamma_{n+1}$ incident to $v_{i_n}$ with 2A and 2B respectively, for some $i$. Each choice of labels 1A, 1B or 2A, 2B indicate a possible choice of "start" or "end" of an MDS sequence.
 %These labels of the ends are given according to the original path orientation of the assembly graph at each vertex. 
To denote a specific HPP with chosen ends path, use $P(\gamma_1,\gamma_{n+1})$. Each $\gamma_1$ and $\gamma_{n+1}$ will be either one of the labels of ending edges and P is composed of complete edges. This path will take upon its own orientation; separate from the orientation of the assembly graph. To obtain the micronuclear sequence, compare the orientation of the HPP to the original path orientation of the assembly graph. These orientations will indicate the MDSs' orientation. Let, $\lambda(e_i)$ be the label for $e_i$ for some $i$, then the micronuclear sequence corresponds to $\lambda(e_o)$,$\ldots$,$\lambda(e_{2n})$. All other edges not visited by HPP are labeled with IESs. Note, an edge may be both 1A and 2B or 1B and 2A (See Figure 7). Such edge is not in the list of HPP edges, $\gamma_1$,$\dots$,$\gamma_{n+1}$. Hence, it also has a label $I_k$ for an IES, for some $k$, so one of the edges may have up to three labels 1A, $I_k$, 2B for example, see Figure 2. In the micronuclear sequence, if the orientation of $\gamma_i$ obtained from the transversal is opposite, the orientation $\gamma_i$ obtains from HPP, label $\gamma_i$ with $\overline{M_{i}}$ otherwise, its label remains $ M_{i}$. We assume that the assembly graphs are open ended, so we will always start and end with IESs in the micronuclear sequences. But, to find all possible micronuclear sequences, we will find all possible combinations of orientations of each HPP with orientations of assembly graph. So, each HPP has four distinct combinations of orientations, with respect to orientation of the HPP and assembly graph. The four distinct orientations are $\Gamma$, $\Gamma^R$, $\Gamma^-$ and $\Gamma^{-R}$. $\Gamma$ is an HPP, $P(\gamma_1,\gamma_{n+1})$. $\Gamma^R$ is the opposite orientation of the assembly graph and same HPP orientation. $\Gamma^-$ is the same orientation as the assembly graph and opposite orientation of HPP. $\Gamma^{-R}$ is both opposite orientations of assembly graph and HPP. To obtain the micronuclear sequences of the other 3 distinct orientations, we go through a similar process of the HPPs' of the following figures, once the orientations of the assembly graph and HPP are fixed. The following HPPs illustrations are $\Gamma$, with its corresponding micronuclear sequence. 
\begin{figure}[htb] 
\subsection{Possible Micronuclear Sequences Represented by Assembly Graph 1212}
\center
\vspace{-8pt}
    \includegraphics[width=10cm] {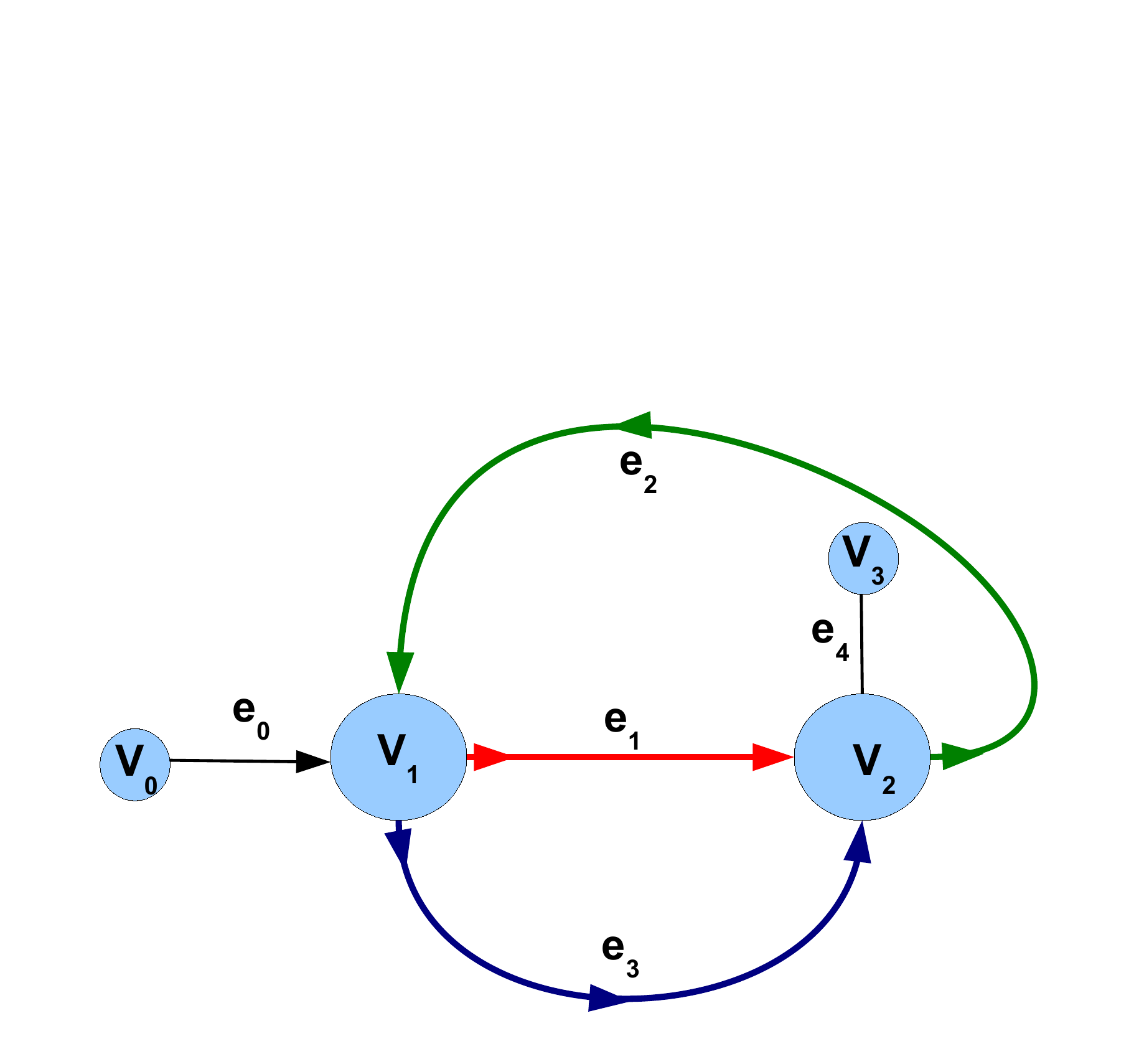}
\vspace{-20pt}
\caption{Assembly graph 1212 consists of the two 4-valent vertices and five edges and transverse path starts with $v_0$,$\dots$,$v_3$. For example, the graph above starts with an edge $e_0$ that first visits $v_1$,$v_2$,$v_1$,$v_2$, giving us assembly graph 1212.}

\label{Assembly1212}
  \end{figure}\\
%The following 2 sequences correspond to Figure~\ref{Assembly1212}.
\newpage
\begin{figure}[h]
\begin{minipage}{0.5\linewidth}
\centering
$1) e_2(1A, 2A)$\\
 $I_0M_1I_1\overline{M_2}I_2 \overline{M_3}I_3$\\
\centering
\includegraphics[width=3in]{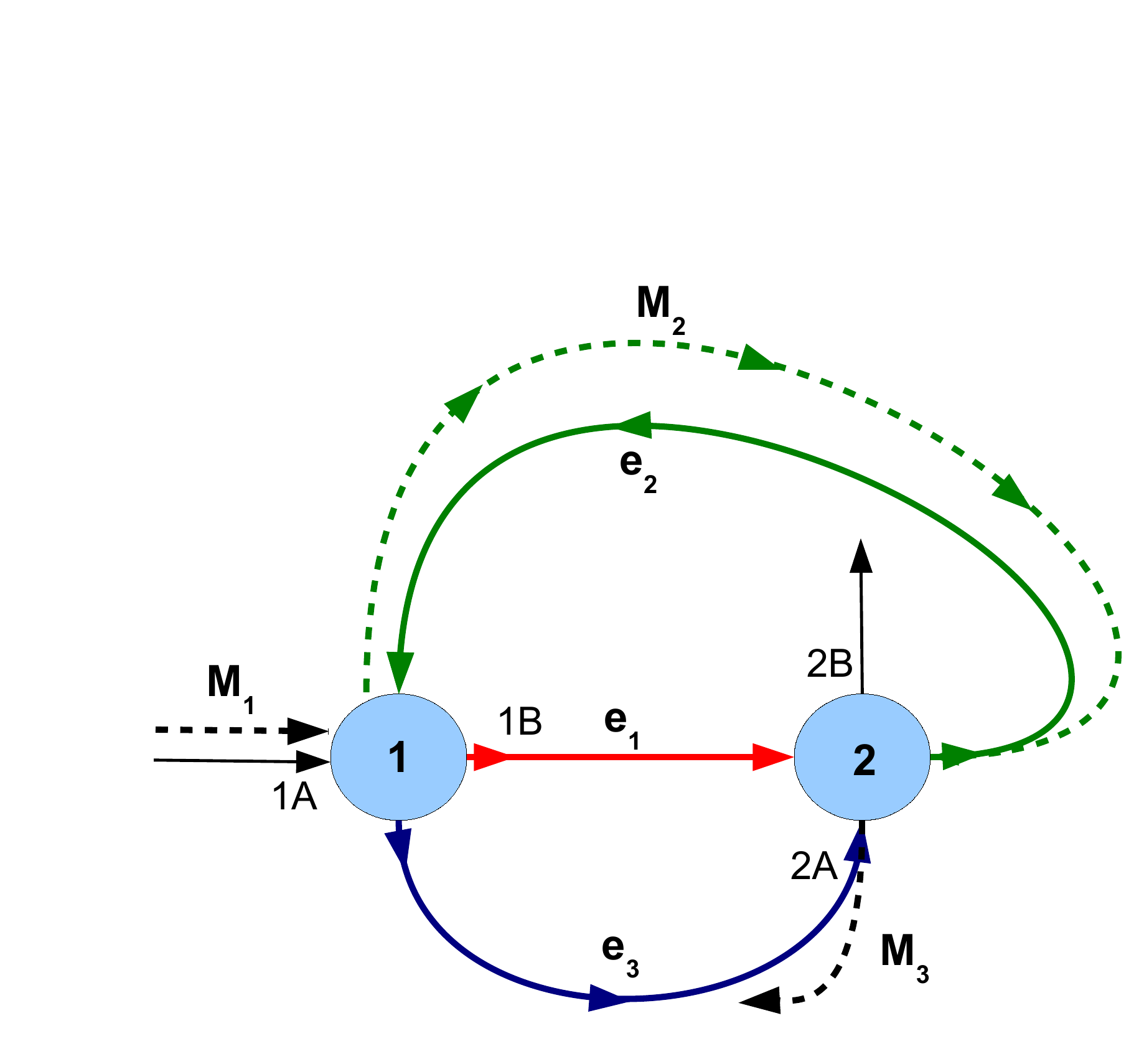}
\centering
Figure 2
\end{minipage}
\hspace{0.5cm}
\begin{minipage}{0.5\linewidth}
\centering
$2 ) e_2(1A, 2B)$\\
$I_0M_1I_1\overline{M_2}I_2 M_3I_3$
\centering
\includegraphics[width=3in]{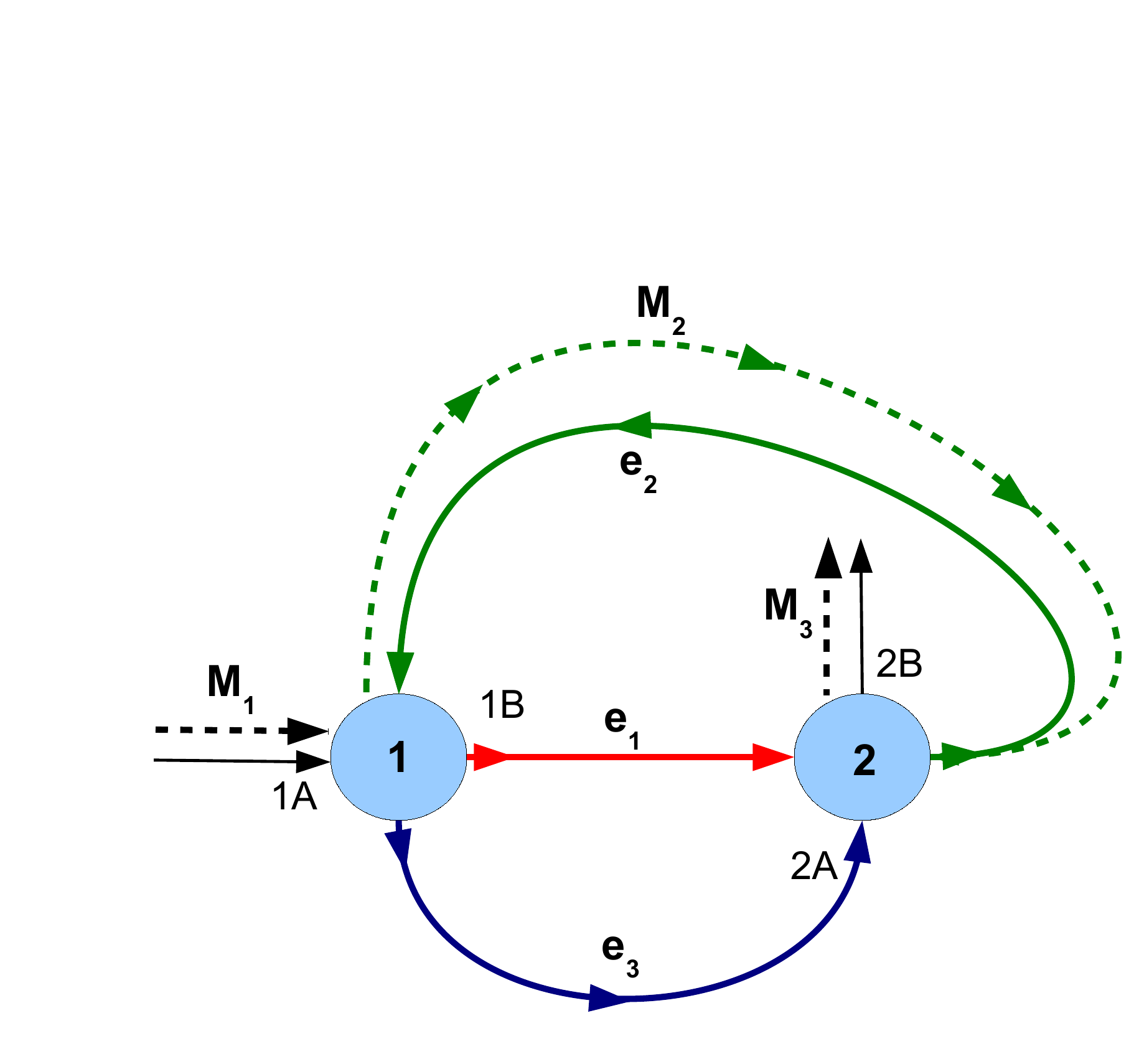}
\centering
Figure 3
\end{minipage}
\end{figure}
\begin{figure}[H]
\begin{minipage}{0.5\linewidth}
\centering
$3) \indent e_2(1B, 2B)$\\
 $I_0\overline{M_1}I_1\overline{M_2}I_2M_3I_3$
\centering
\includegraphics[width=3in]{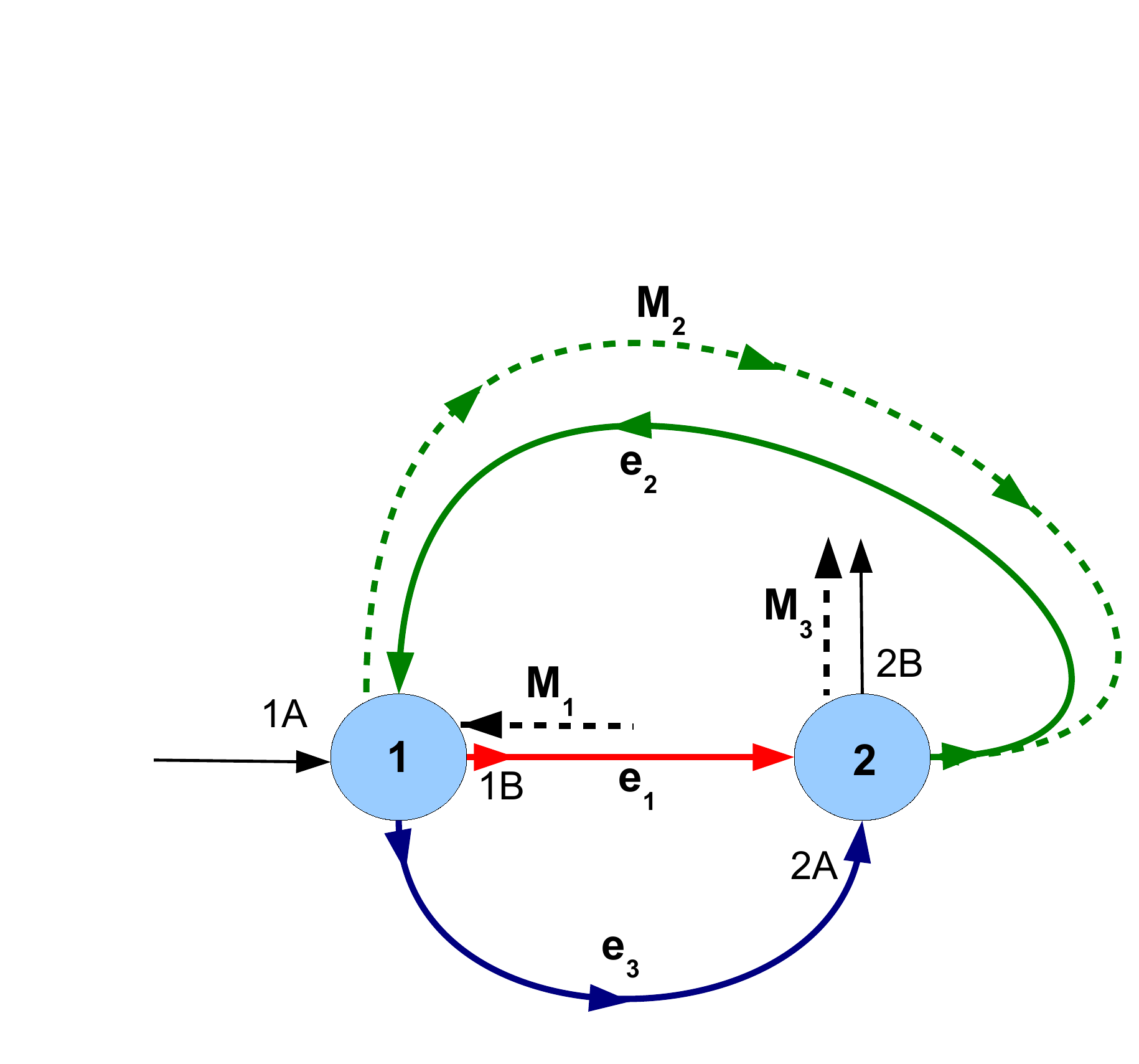}
\centering
Figure 4
\end{minipage}
\hspace{0.5cm}
\begin{minipage}{0.5\linewidth}
\centering
$4) \indent e_2(1B, 2A)$\\
$I_0\overline{M_1}I_1\overline{M_2}I_2\overline {M_3}I_3$\\
\centering
\includegraphics[width=3in]{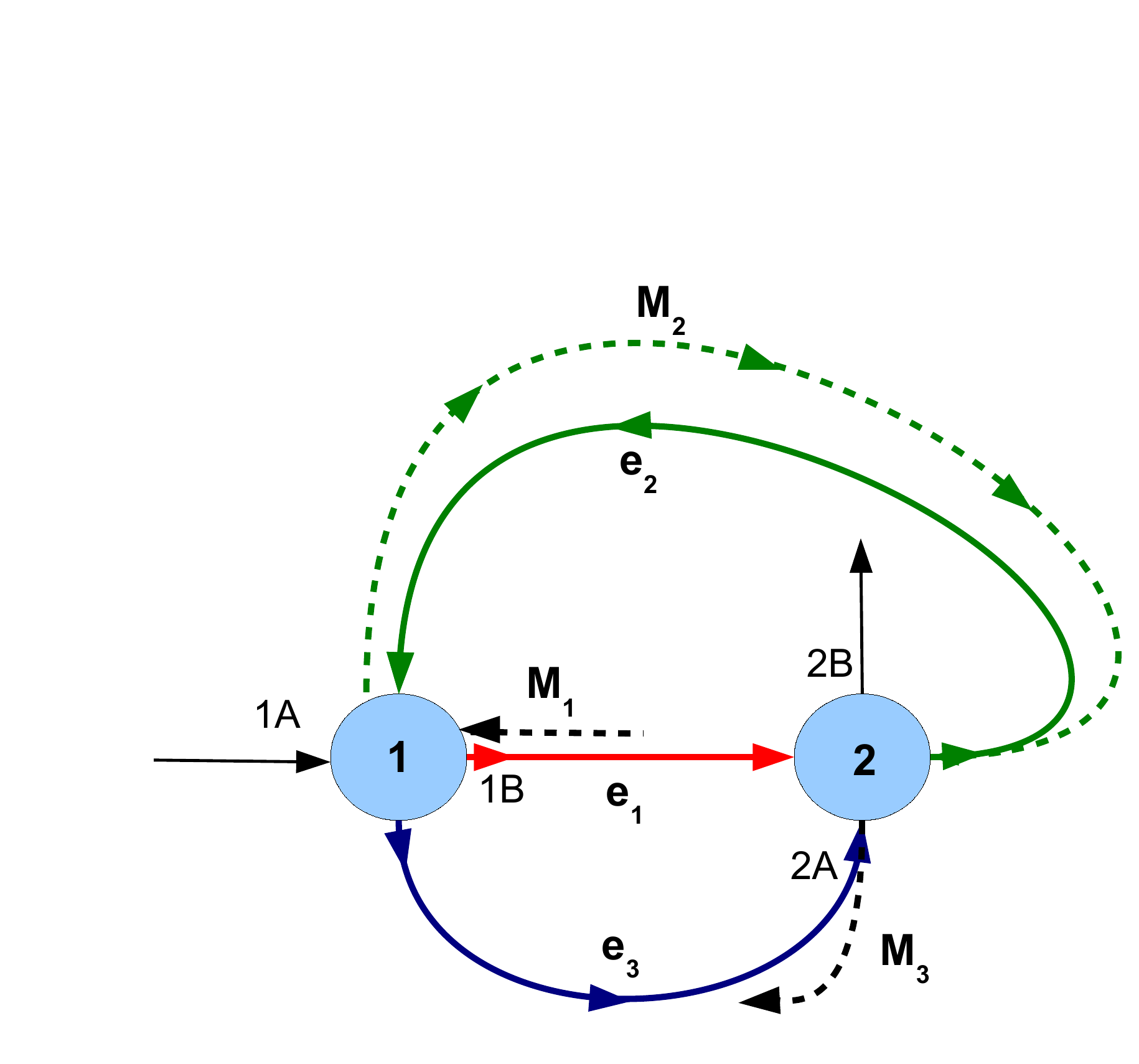}
\centering
Figure 5
\end{minipage}
\end{figure}
\begin{center}
The figures above are HPPs ($\Gamma$) grouped consisting of the path ($e_2$).\\
\end{center}
Figure 2 has two vertices so there is 1 full edge in the HPP ($e_2$), with ending edges ($1A, 2A$).  Full edge, in the sense that it is a complete edge of the assembly graph and $\Gamma$. So, the first line of notation above Figure 2 denotes HPP as $ e_2(1A, 2A)$. The second line of notation above the assembly graph represents the micronuclear sequence of the Hamiltonian polygonal path. To obtain the micronuclear sequence, read the labels of the transversal  path of the given graph. For example, we always start with an IES $(I_0)$ following $1A$, which is the first MDS ($M_1$) and continues along $e_1$, which is $I_1$. Now, moving along $e_2$  $(\overline{M_2})$ we continue to $I_2$ which is part of $e_3$. Then, with $\overline{M_3}$ ($2A$), however, we must complete the assembly graph and micronuclear sequence with an IES ($I_3$). Notice, the MDS sequences in the HPP are inverted, because their orientations are the opposite with respect to the transversal. Obtaining the micronuclear sequences for $\Gamma^R$, $\Gamma^-$ and $\Gamma^{-R}$ is the same process as obtaining $\Gamma$, as long as the combinations of orientations for each are distinct, with respect to the assembly graph and HPP. \\
\begin{figure}[h]
\begin{minipage}{0.5\linewidth}
\vspace{-20pt}
\centering
$5) \indent e_3(1A, 2A)$\\
 $I_0M_1I_1\overline{M_3}I_2 M_2I_3$\\
\centering
\includegraphics[width=3in]{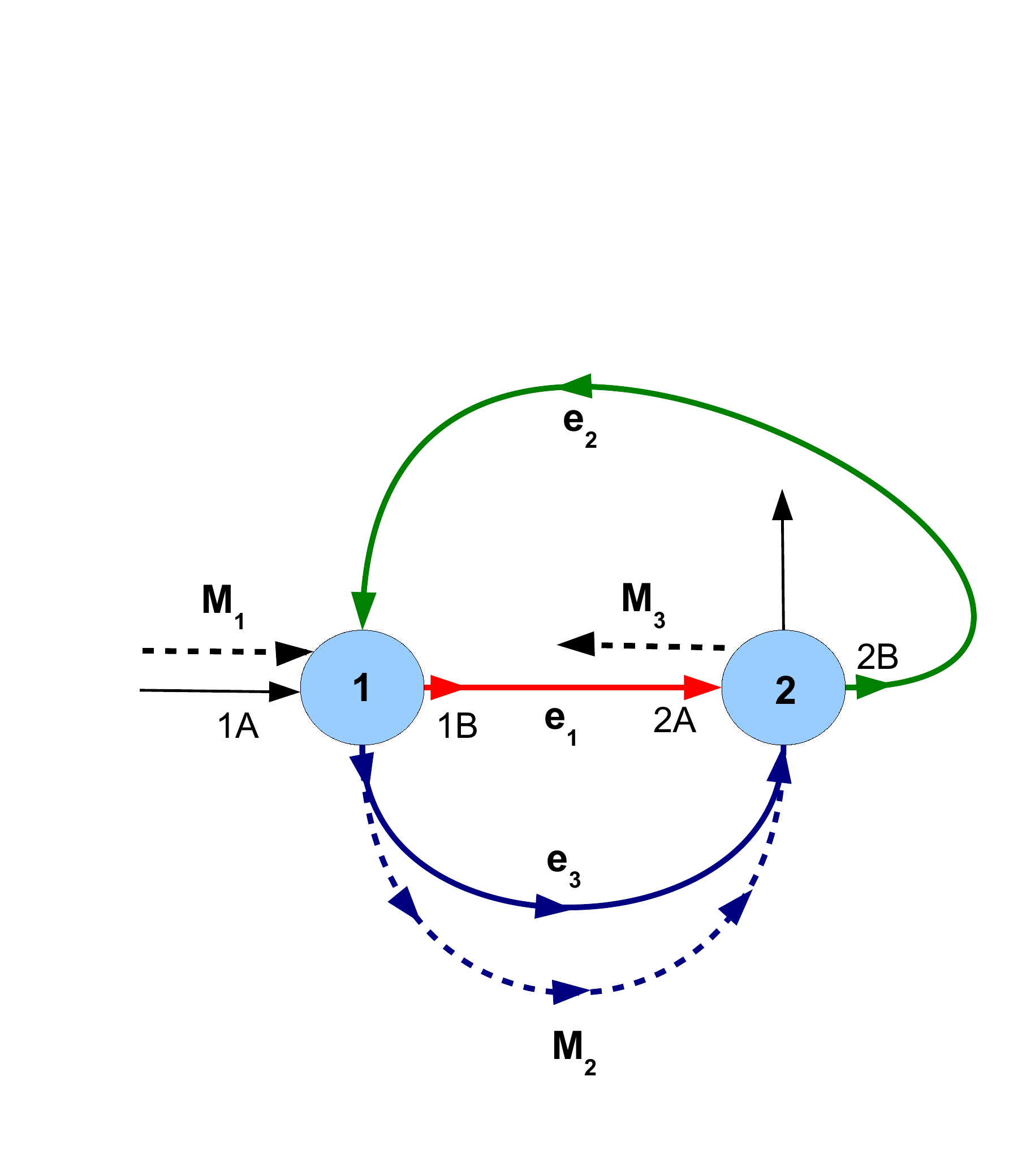}
\vspace{10pt}
\centering
Figure 6
\end{minipage}
\hspace{0.5cm}
\begin{minipage}{0.5\linewidth}
\vspace{-20pt}
\centering
$6) \indent e_3(1A, 2B)$\\
$I_0M_1I_1 M_3 I_2 M_2I_3$\\
\centering
\includegraphics[width=3in]{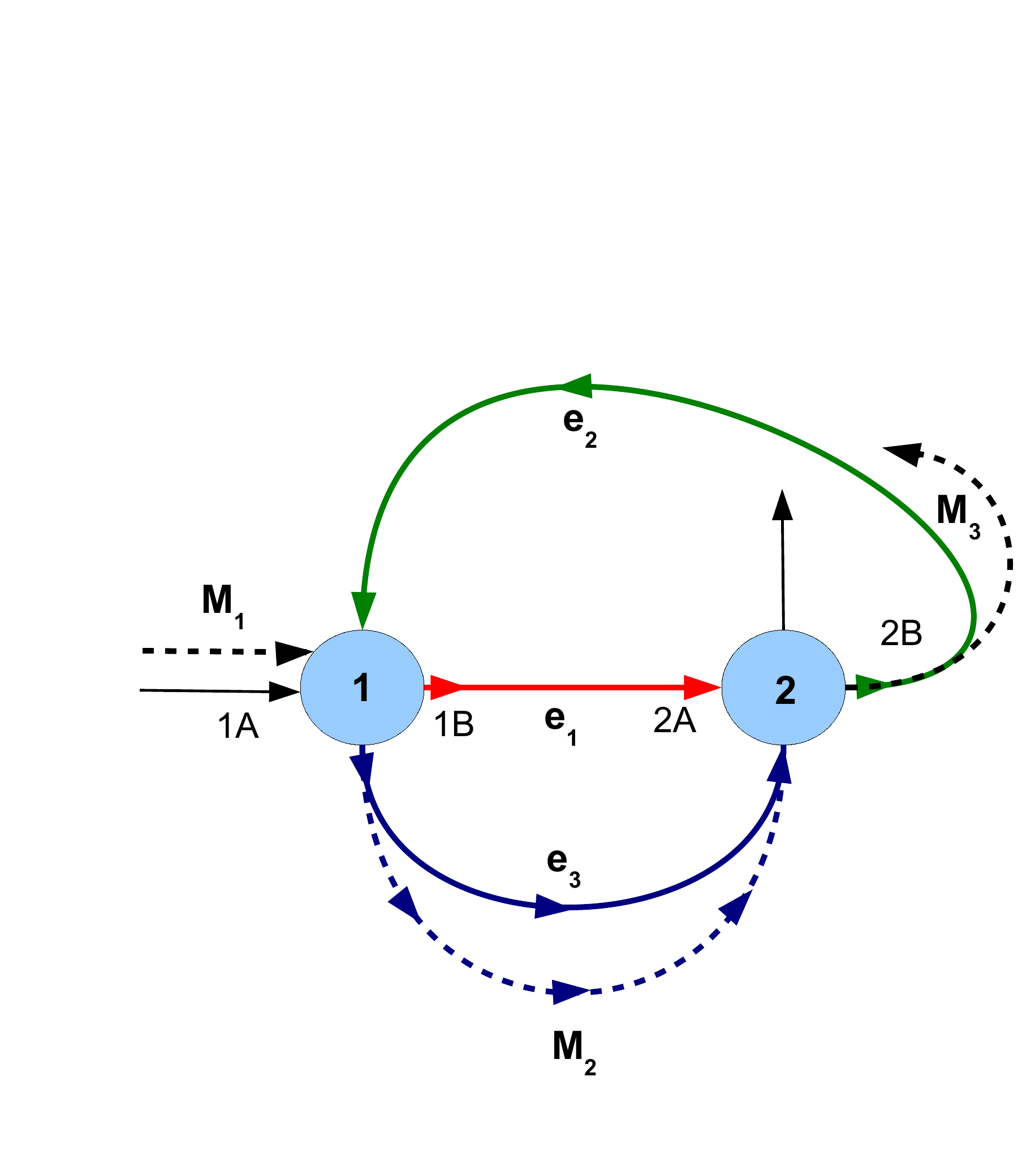}
\vspace{10pt}
\centering
Figure 7
\end{minipage}
\end{figure}
\begin{figure}[H]
\begin{minipage}{0.5\linewidth}
\vspace{-33pt}
\centering
$7) \indent e_3(1B, 2A)$ \\
$I_0\overline{M_1}I_1\overline{M_3}I_2 M_2 I_3$\\
\centering
\includegraphics[width=3in]{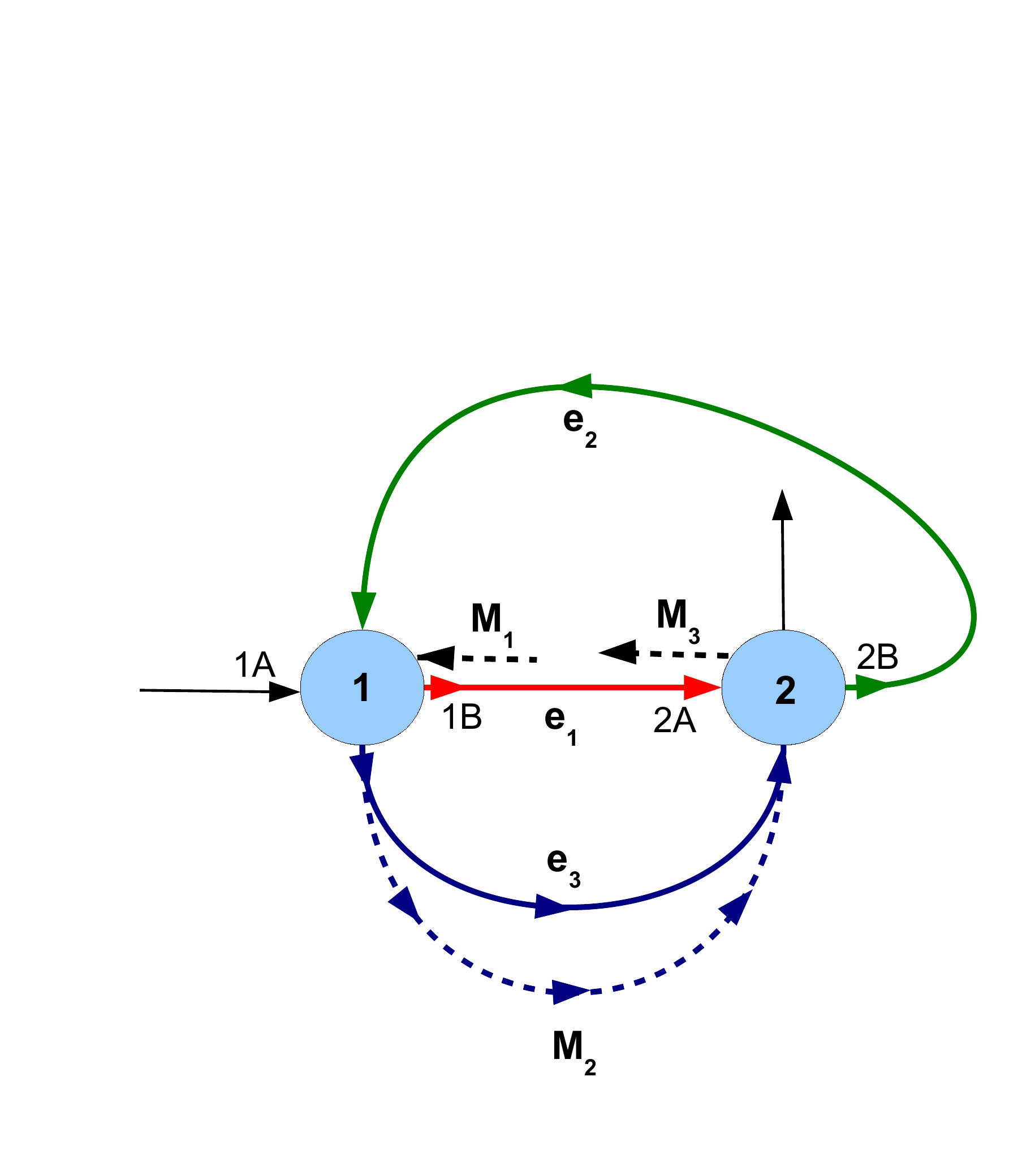}
\centering
Figure 8
\end{minipage}
\hspace{0.5cm}
\begin{minipage}{0.5\linewidth}
\vspace{-33pt}
\centering
$8) \indent e_3(1B, 2B)$\\
$I_0\overline{M_1}I_1 M_3 I_2 M_2I_3$\\
\centering
\includegraphics[width=3in]{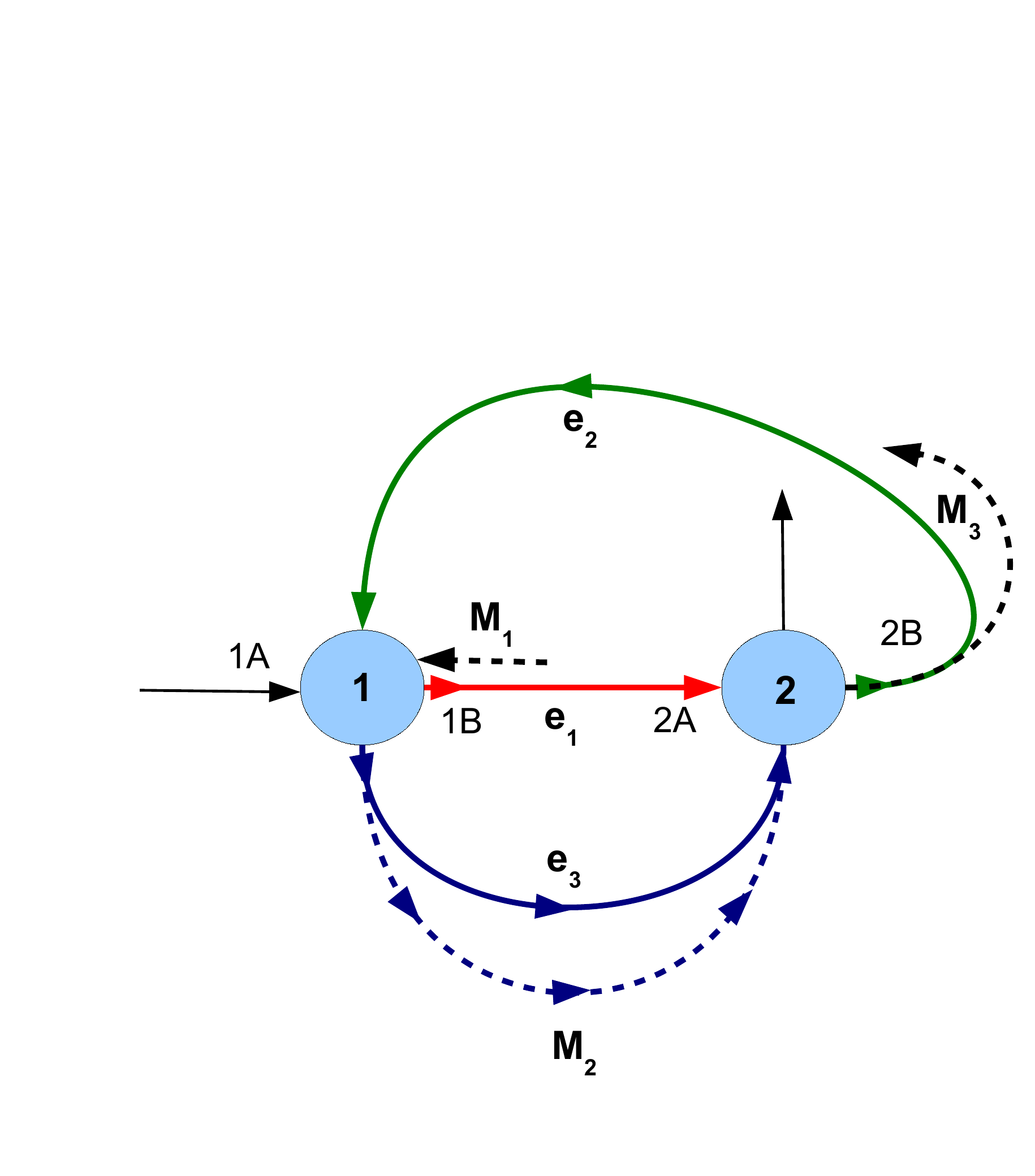}
\centering
Figure 9
\end{minipage}
\end{figure}
\begin{center}
The figures above are HPPs ($\Gamma$) grouped consisting of the path ($e_3$).\\
\end{center}
\newpage
\begin{figure}[h]
\begin{minipage}{0.5\linewidth}
\centering
$9) \indent e_1(1A, 2A)$\\
\indent \hspace{6mm} $I_0 M_2 I_1 M_1 I_2 \overline{M_3}I_3$\\
\centering
\includegraphics[width=3in]{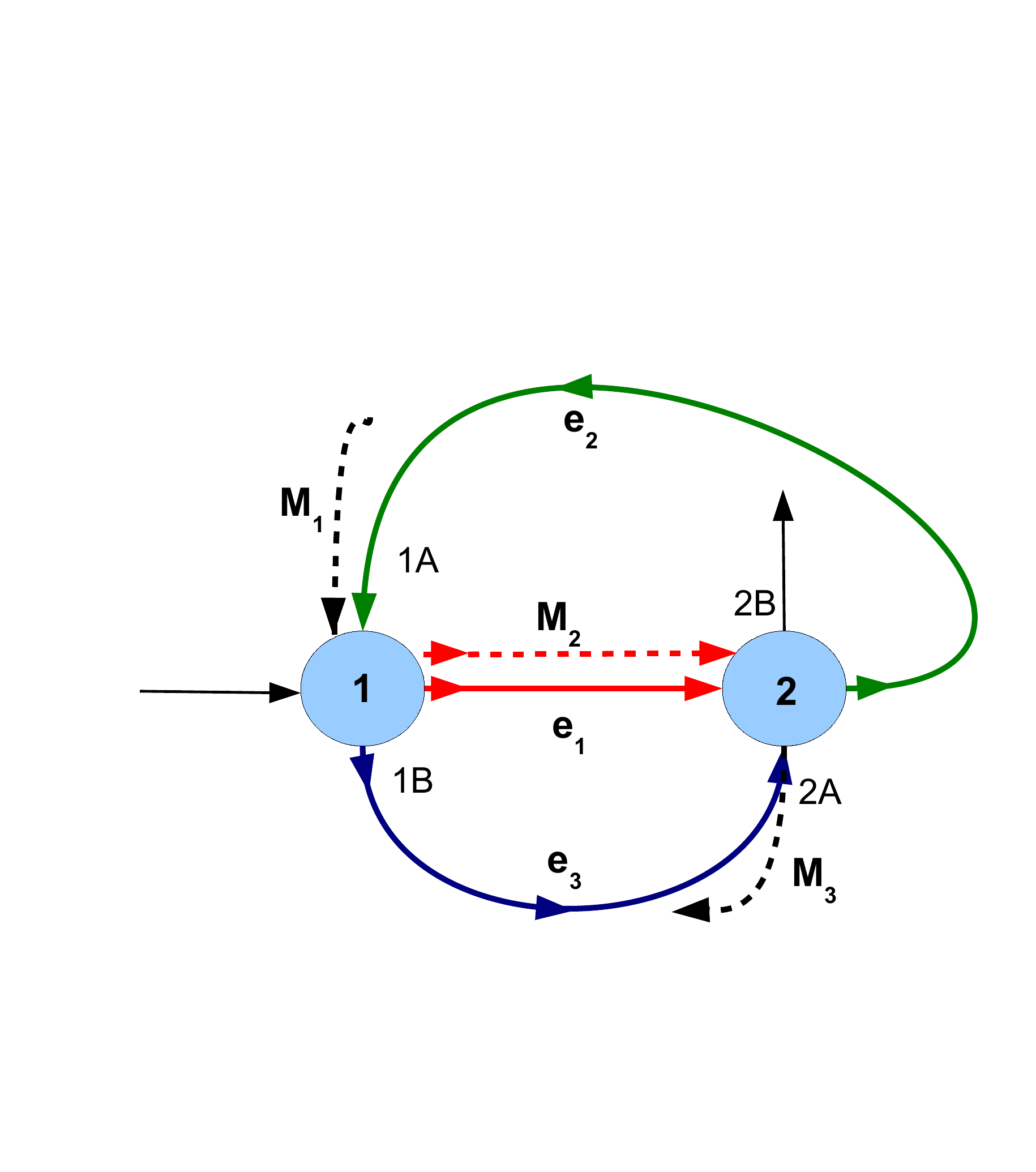}
\centering
Figure 10
\end{minipage}
\hspace{0.5cm}
\begin{minipage}{0.5\linewidth}
\centering
$10) \indent e_1(1A, 2B)$\\
 $I_0 M_2 I_1 {M_1} I_2 {M_3}I_3$\\
\centering
\includegraphics[width=3in]{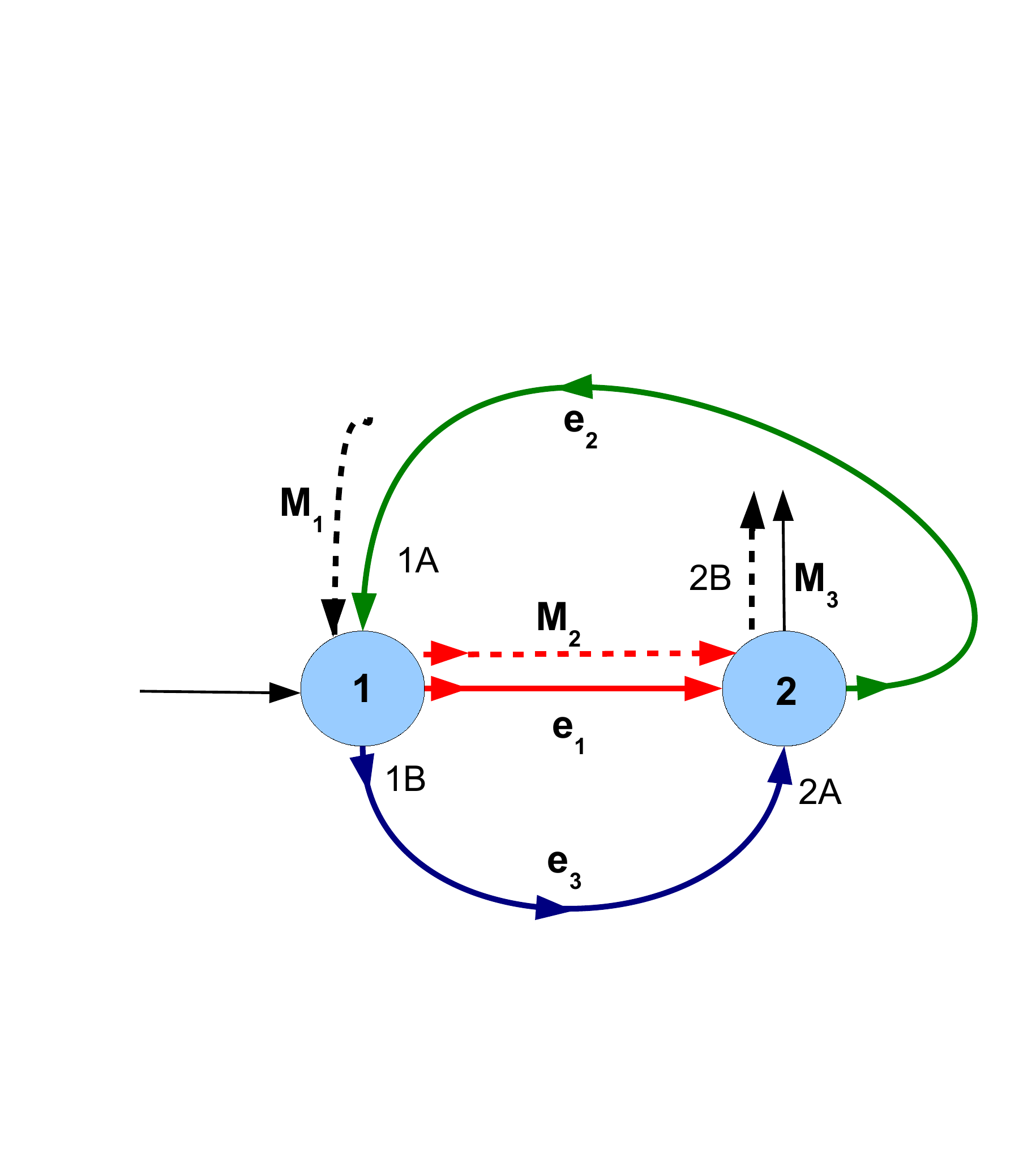}
\centering
Figure 11
\end{minipage}
\end{figure}
\begin{figure}[H]
\begin{minipage}{0.5\linewidth}
\centering
$11) \indent e_1(1B, 2A)$\\
\indent \hspace{6mm} $I_0 M_2 I_1\overline M_1 I_2 \overline{M_3}I_3$\\
\centering
\includegraphics[width=3in]{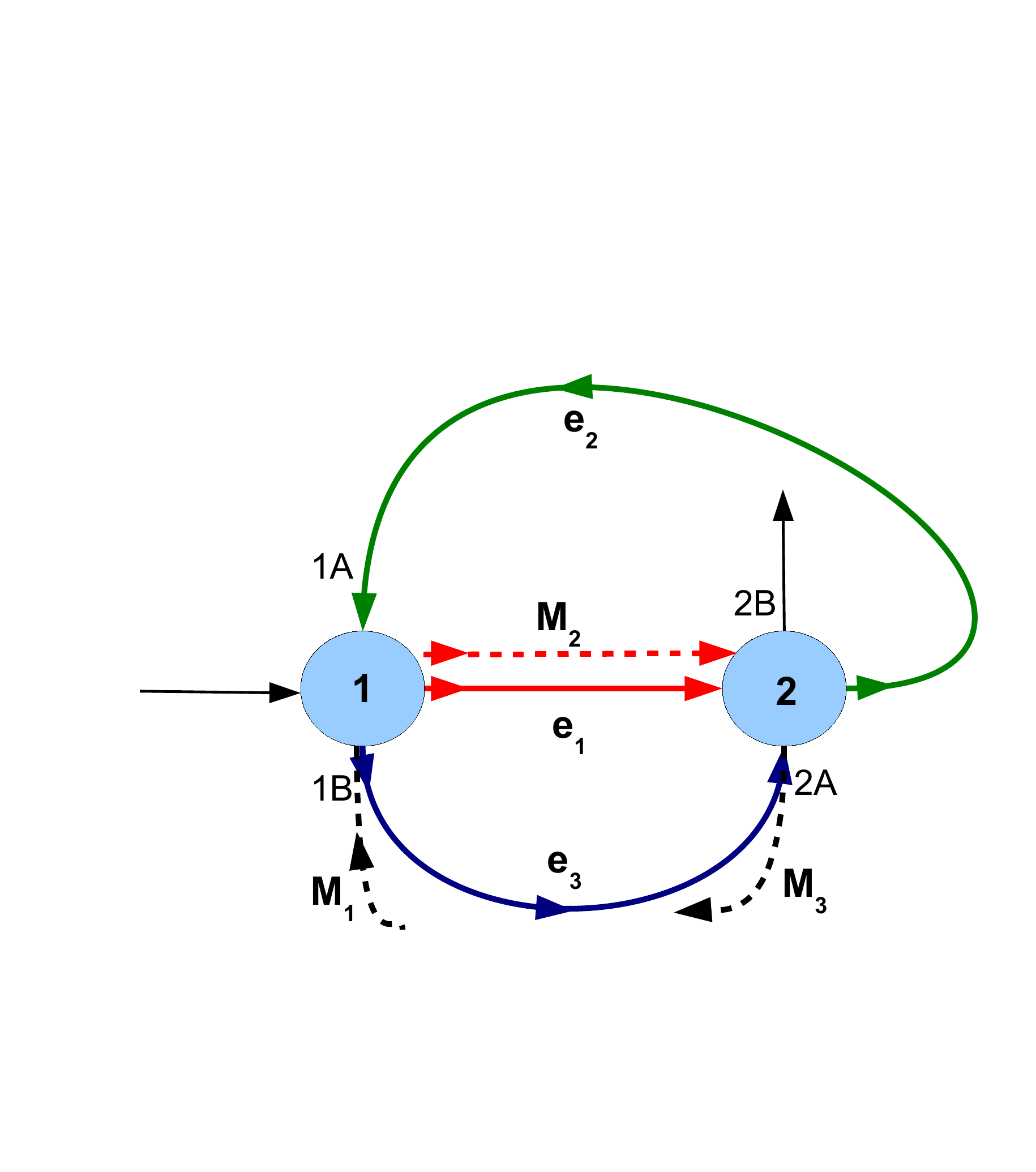}
\centering
Figure 12
\end{minipage}
\hspace{0.5cm}
\begin{minipage}{0.5\linewidth}
\centering
$12) \indent e_1(1B, 2B)$\\
$I_0 M_2 I_1\overline M_1 I_2 M_3$
\centering
\includegraphics[width=3in]{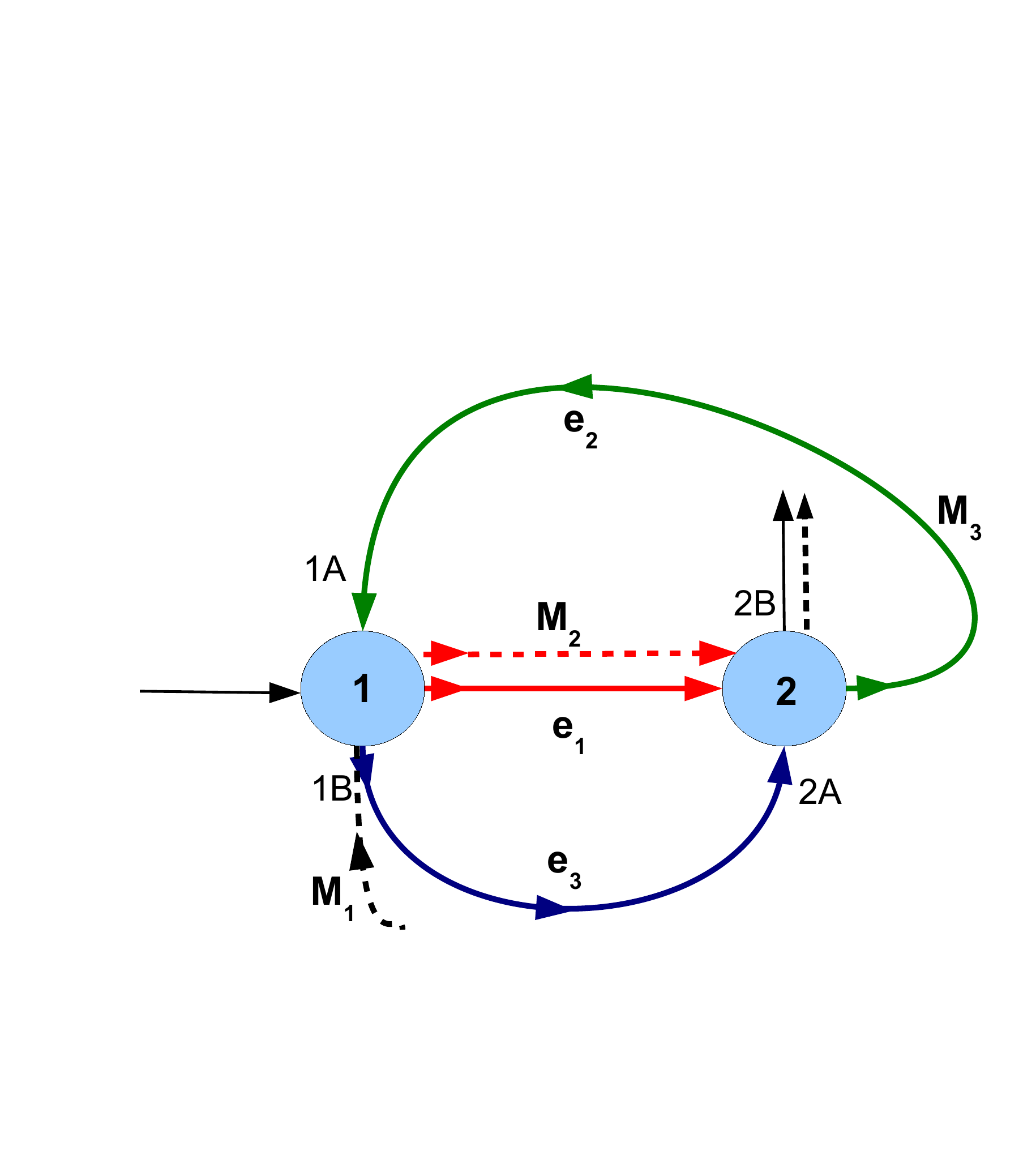}
\centering
Figure 13
\end{minipage}
\end{figure}
\begin{center}
The figures above are HPPs ($\Gamma$)grouped consisting of the path ($e_1$).\\
\end{center}

\begin{table}
\begin{tabular}{|c|c|c|}
\hline
HPP & Micronuclear Sequence & Inverted Micronuclear Sequence\\
\hline
\hline
\multirow{4}{*}{}$1) \indent e_2(1A, 2A)$     & $I_0M_1I_1\overline{M_2}I_2 \overline{M_3}I_3$ &$I_0 M_3I_1M_2 I_2\overline{M_1}I_3$\\
$2 )\indent  e_2(1A, 2B)$ & $I_0M_1I_1\overline{M_2}I_2 M_3I_3$ & $I_0\overline{M_3}  I_1M_2I_2 \overline{M_1}I_3$ \\
$3) \indent e_2(1B, 2B)$ & $I_0\overline{M_1}I_1 \overline{M_2}I_2M_3I_3$ &$I_0\overline{M_3}  I_1 M_2 I_2 M_1 I_3$\\
$4) \indent e_2(1B, 2A)$ & $I_0\overline{M_1}I_1\overline{M_2}I_2\overline {M_3}I_3$ & $I_0 M_3 I_1 M_2 I_2 M_1I_3$\\
\hline
\multirow{4}{*}{}$5)\indent e_3(1A, 2A)$ & $I_0M_1I_1\overline{M_3}I_2 M_2I_3$ & $I_0 \overline{M_2} I_1 M_3 I_2 \overline{M_1}I_3$\\
$6) \indent e_3(1A, 2B)$ &$I_0M_1I_1 M_3 I_2 M_2I_3$ & $I_0\overline{M_2}I_1\overline{M_3}I_2\overline{M_1}I_3$\\
 $7) \indent e_3(1B, 2A)$ & $I_0\overline{M_1}I_1\overline{M_3}I_2 M_2 I_3$ & $I_0\overline{M_2}I_1 M_3 I_2 M_1 I_3$\\
$8) \indent e_3(1B, 2B)$ & $I_0\overline{M_1}I_1 M_3 I_2 {M_2}I_3$ & $I_0\overline{M_2} I_1 \overline{M_3} I_2 M_1 I_3$\\
\hline
$9) \indent e_1(1A, 2A)$ & $I_0 M_2 I_1 M_1 I_2 \overline{M_3}I_3$ & $I_0 M_3 I_1\overline{M_1} I_2 \overline{M_2} I_3$\\
$10) \indent e_1(1A, 2B)$  & $I_0 M_2 I_1 {M_1} I_2 {M_3}I_3$ &  $ I_0 \overline{M_3} I_1\overline{M_1} I_2 \overline{M_2}I_3$\\
$11) \indent  e_1(1B \hspace{2mm}2A)$ & $I_0 M_2 I_1\overline M_1 I_2\overline{M_3}I_3$ & $I_0 M_3I_1M_1I_2\overline{M_2} I_3$\\
$12) \indent  e_1(1B \hspace{2mm}2B)$ & $I_0 M_2 I_1\overline M_1 I_2 M_3I_3$ & $I_0\overline{M_3}I_1 M_1 I_2 \overline{M_2} I_3$ \\
\hline
\end{tabular}
\end{table}
\newpage
The first column in the table above are the possible HPPs ( $\Gamma$) of assembly graph 1212. The second column represent the micronuclear sequences obtained from the HPPs. The third column is $\Gamma^R$. Some sequences in $\Gamma^-$ and $\Gamma^{-R}$ corresponded to $\Gamma$, $\Gamma^R$, which leaves us with twenty four distinct labellings of the assembly graph 1212. Notice, that $\Gamma^R$ is complementary to the original, $\Gamma$  and therefore represent the same DNA sequence \cite{2}. Graphically, the difference between  $\Gamma^R$ and $\Gamma$  is the orientation of the assembly graph, which inverts the sequences when obtained.
\newpage
\subsection{Possible Micronuclear Sequences Represented by Assembly Graph 1221}
\begin{figure}[h]
\centering
\includegraphics[width= 10cm]{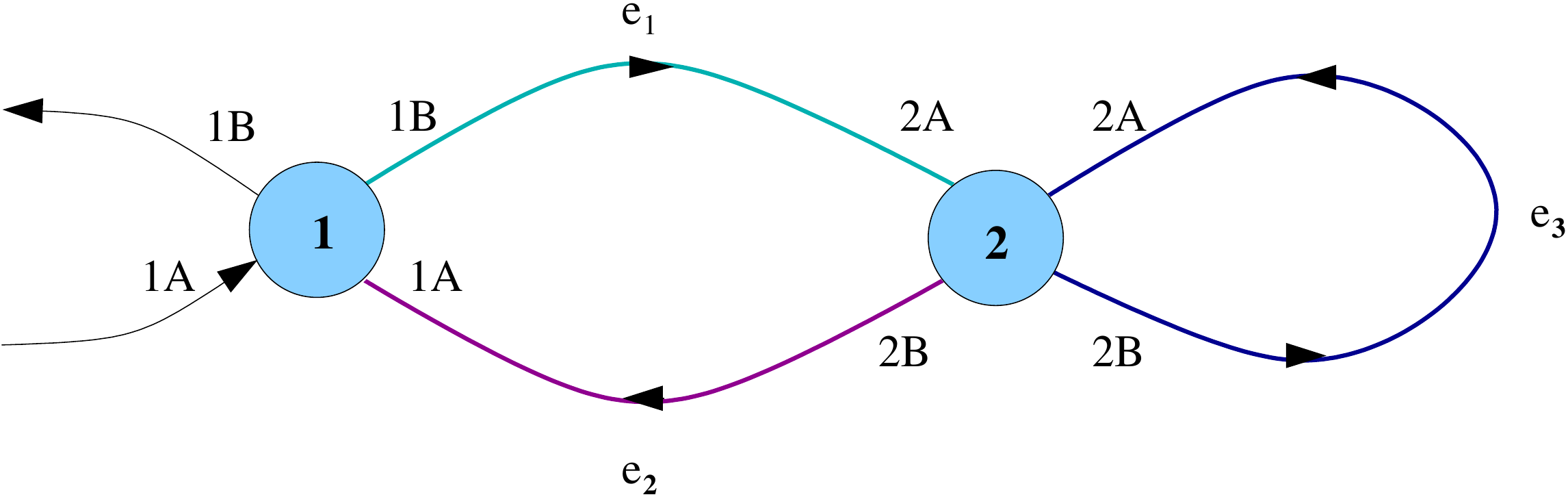}
\caption{The assembly graph above is 1221, which is formed similarly to the previous assembly graph.} 
\end{figure}
\begin{figure}[h]
\begin{minipage}{0.5\linewidth}
 $1) \indent e_1(1A, 2A)$\\
$I_0M_2 I_1 \overline{M_3}I_2 M_1I_3$\\

\centering
\includegraphics[width=3in]{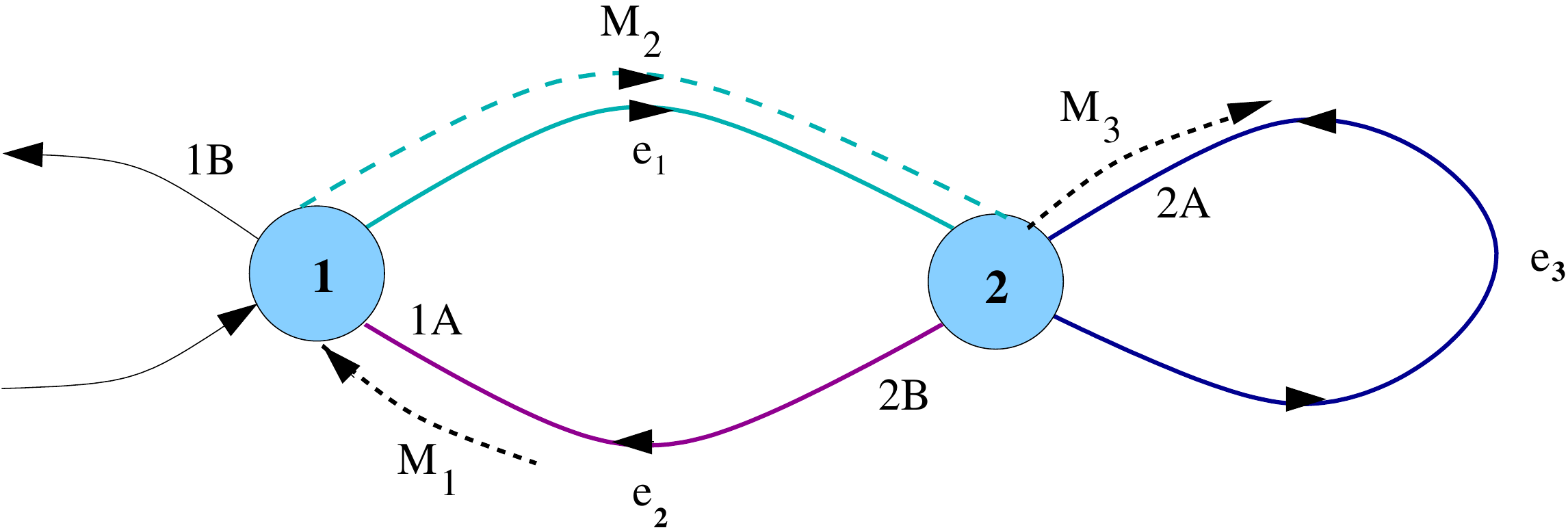}
\centering
Figure 2
\end{minipage}
\hspace{0.5cm}
\begin{minipage}{0.5\linewidth}
$2) \indent e_1(1A, 2B)$\\
$I_0M_2I_1 M_3I_2M_1I_3$

\centering
\includegraphics[width=3in]{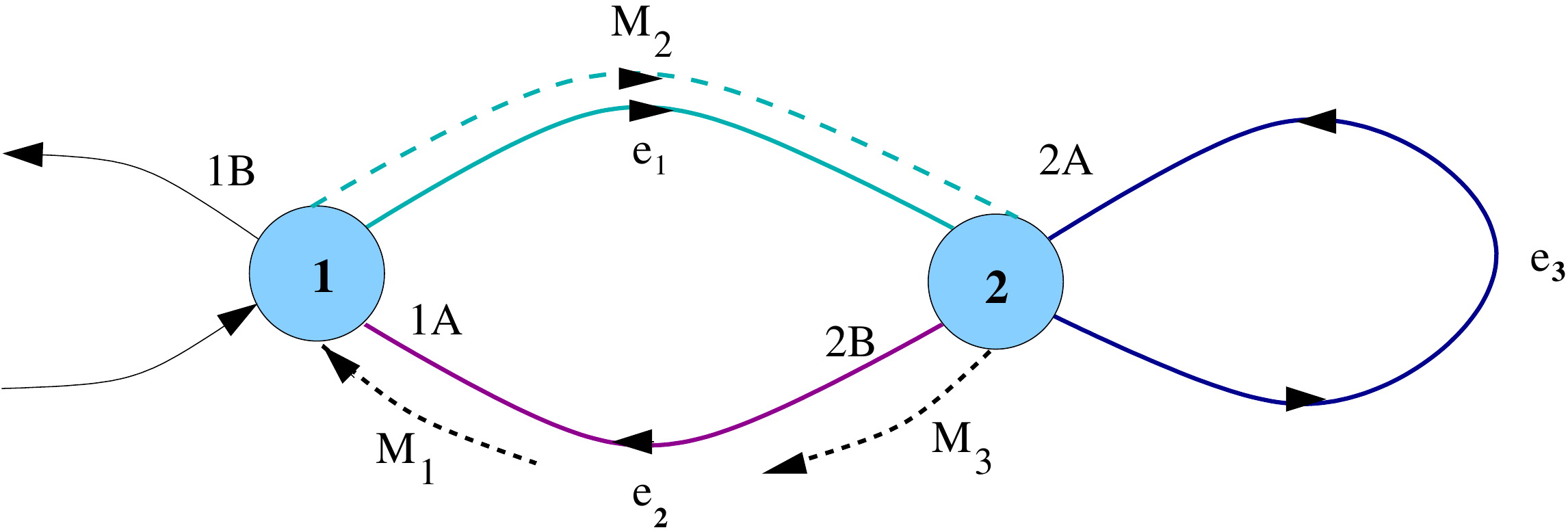}
\centering
Figure 3
\end{minipage}
\end{figure}
\begin{figure}[H]
\begin{minipage}{0.5\linewidth}
$3) \indent e_1(1B, 2A)$\\
$I_0M_2I_1\overline{ M_3}I_2\overline{M_1}I_3$\\

\centering
\includegraphics[width=3in]{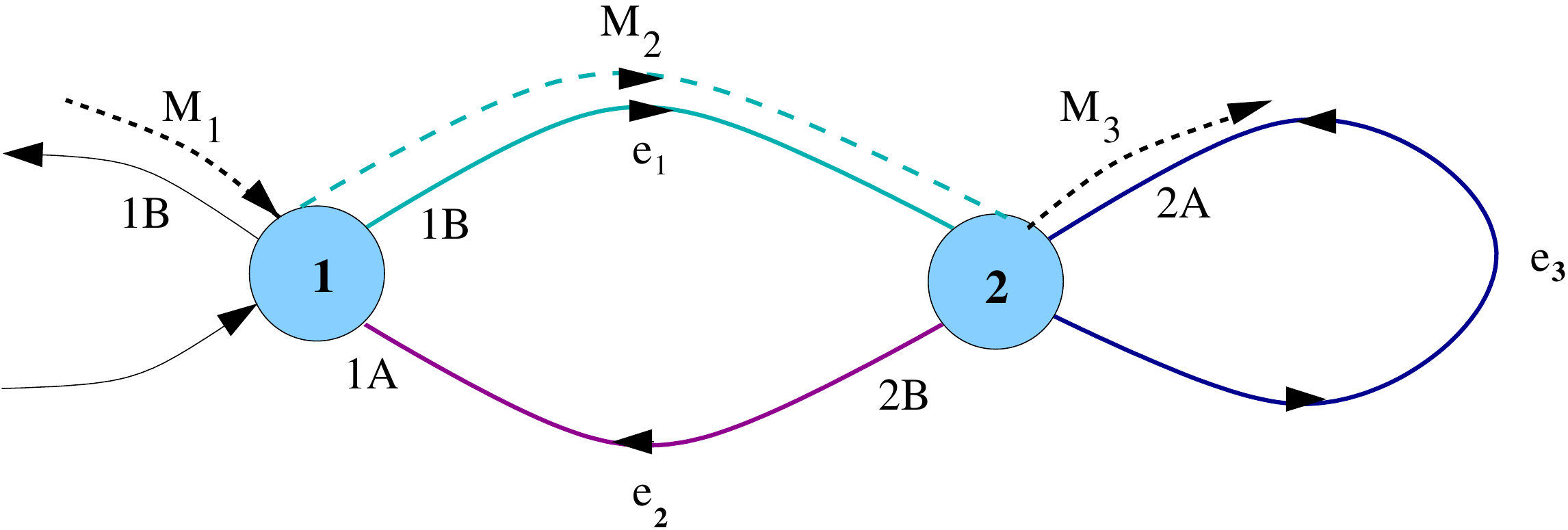}
\centering
Figure 4
\end{minipage}
\hspace{0.5cm}
\begin{minipage}{0.5\linewidth}
$4) \indent e_1(1B, 2B)$\\
$I_0 M_2 I_1M_3 I_2 \overline{M_1}I_3$\\

\centering
\includegraphics[width=3in]{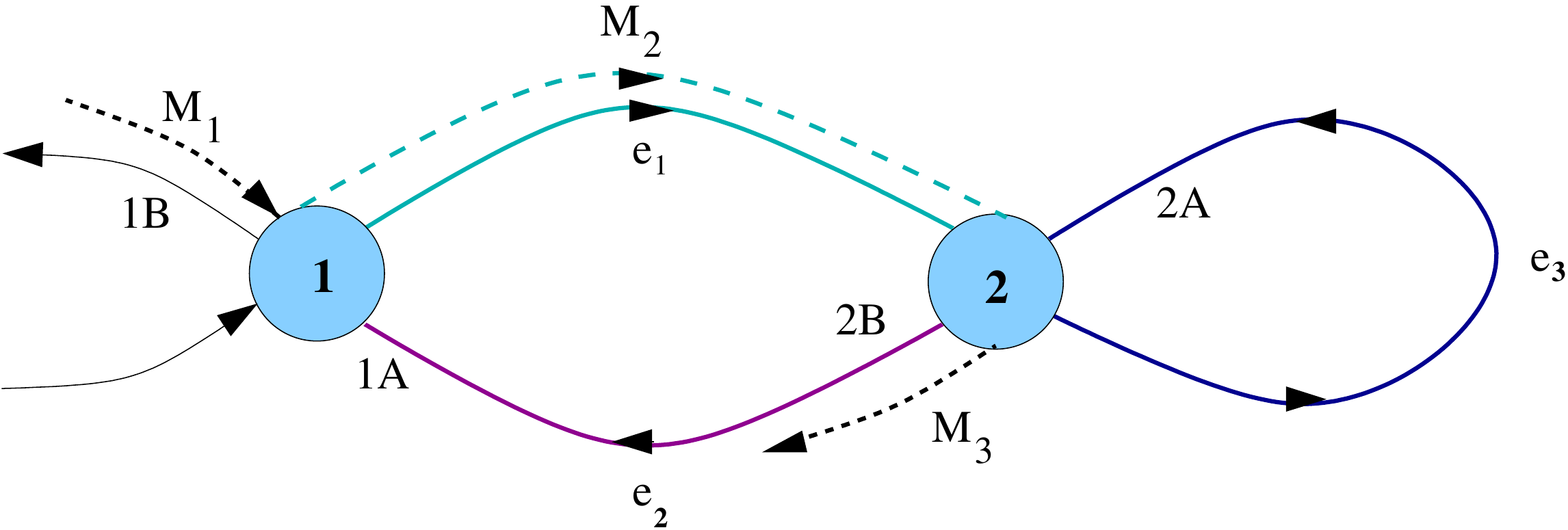}
\centering
Figure 5
\end{minipage}
\end{figure}
The figures above are HPPs ($\Gamma$) grouped consisting of the path ($e_1$). The corresponding micronuclear sequences are obtained in a similar fashion as the previous micronuclear sequences of previous assembly graph 1212.\\
\newpage
\begin{figure}[h]
\begin{minipage}{0.5\linewidth}
 $5) \indent e_2(1A, 2B)$\\
$I_0M_1I_1M_3 I_2\overline{M_2} I_3$\\

\centering
\includegraphics[width=3in]{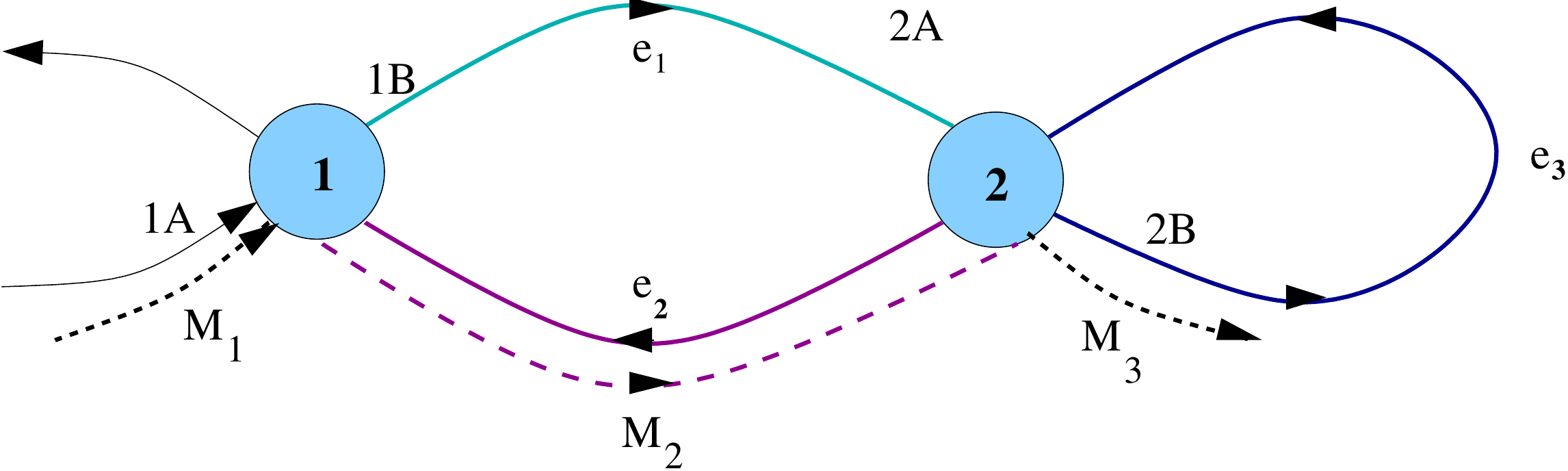}
\centering
Figure 6
\end{minipage}
\hspace{0.5cm}
\begin{minipage}{0.5\linewidth}
$6) \indent e_2(1A, 2A)$\\
$I_0M_1I_1\overline{ M_3} I_2\overline{M_2}I_3$\\

\centering
\includegraphics[width=3in]{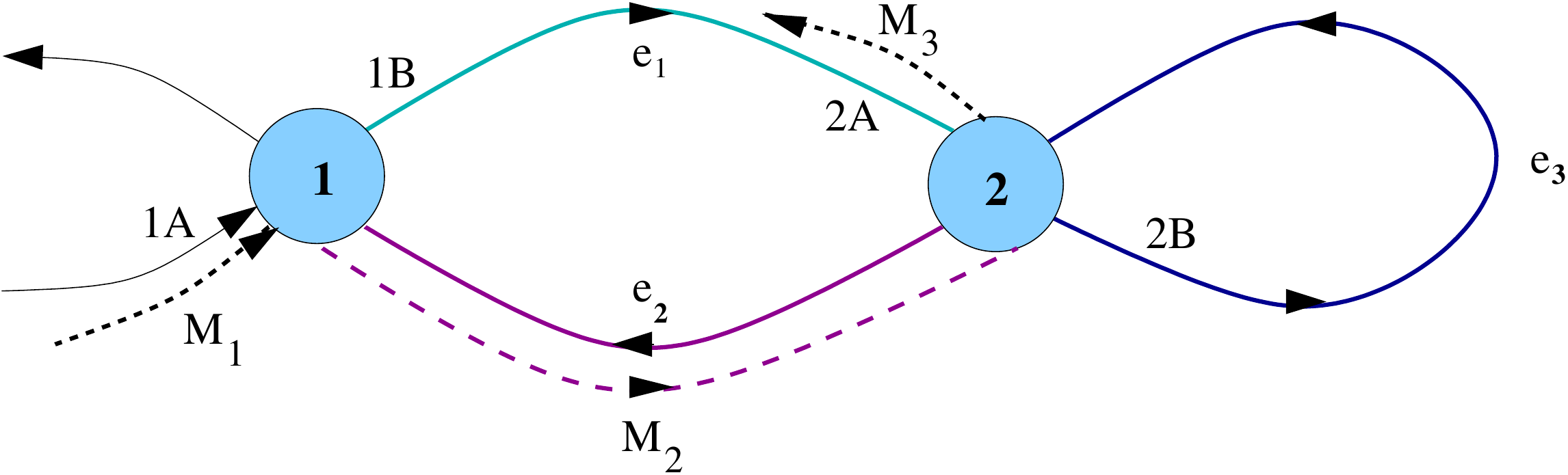}
\centering
Figure 7
\end{minipage}
\end{figure}
\begin{figure}[H]
\begin{minipage}{0.5\linewidth}
$7) \indent e_2(1B, 2A)$\\
$I_0\overline{M_1}I_1\overline{M_3}I_2 \overline{M_2} I_3$\\

\centering
\includegraphics[width=3in]{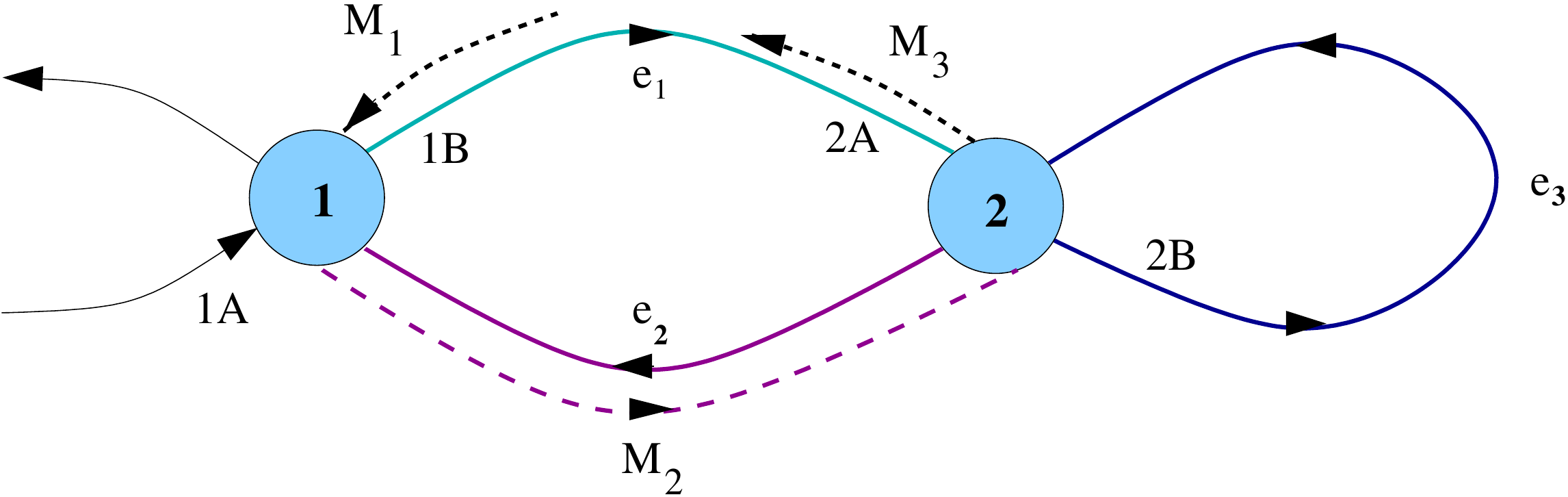}
\centering
Figure 8
\end{minipage}
\hspace{0.5cm}
\begin{minipage}{0.5\linewidth}
$8) \indent e_2(1B, 2B)$\\
$I_0\overline{M_1}I_1 M_3 I_2 \overline{M_2}I_3$\\

\centering
\includegraphics[width=3in]{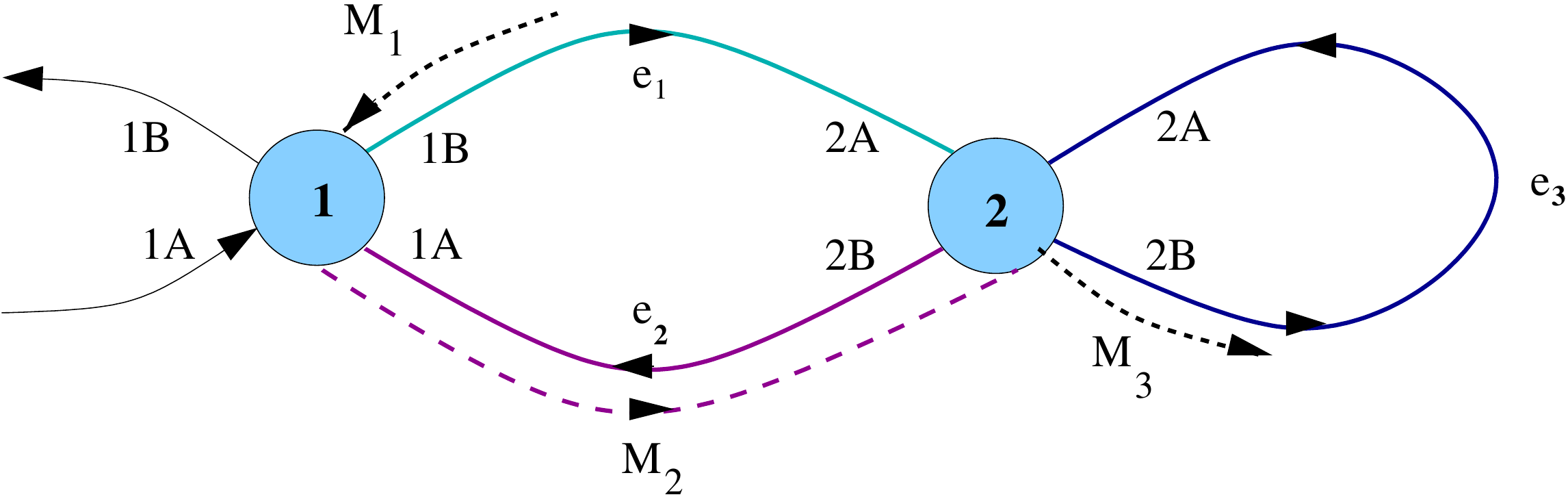}
\centering
Figure 9
\end{minipage}
\end{figure}
The figures above are HPPs ($\Gamma$) grouped consisting of the path ($e_2$). The corresponding micronuclear sequences are obtained in a similar fashion as the previous micronuclear sequences of previous assembly graph 1212.\\
\newpage

\begin{table}
\begin{tabular}{|c|c|c|}
\hline
HPP & Micronuclear Sequence & Inverted Micronuclear Sequence\\
\hline
\hline
\multirow{4}{*}{}$1) \indent e_1(1A, 2A)$ &  $I_0M_2 I_1 \overline{M_3}I_2 M_1I_3$ & $I_0\overline{M_1} I_1M_3I_2 \overline{M_2}I_3$ \\
$2) \indent e_1(1A, 2B)$ & $I_0M_2I_1 M_3I_2M_1I_3$  & $I_0 \overline{M_1}I_1 \overline{M_3}I_2 \overline{M_2}I_3$\\
$3) \indent e_1(1B, 2A)$ & $I_0M_2I_1\overline{ M_3}I_2\overline{M_1}$ & $M_1  I_0 M_3 I_1\overline{ M_2} I_2$\\
$4) \indent e_1(1B, 2B)$ & $I_0 M_2 I_1M_3 I_2 \overline{M_1}$ & $ M_1 I_0 \overline{M_3} I_1\overline{ M_2}I_2$\\
\hline
\multirow{4}{*}{}$5) \indent e_2(1A, 2B)$ &  $I_0M_1I_1M_3 I_2\overline{M_2} I_3$ & $I_0 M_2 I_1 \overline{M_3} I_2 \overline{M_1}I_3$\\
$6) \indent e_2(1A, 2A)$ & $I_0M_1I_1\overline{ M_3} I_2\overline{M_2}I_3$ & $I_0 M_2 I_1 M_3 I_2\overline{M_1}I_3$\\
 $7) \indent e_2(1B, 2A)$ & $I_0\overline{M_1}I_1\overline{M_3}I_2 \overline{M_2} I_3$ & $I_0M_2 I_1 M_3 I_2 M_1 I_3$\\
$8) \indent e_2(1B, 2B)$ & $I_0\overline{M_1}I_1 M_3 I_2 \overline{M_2}I_3$ & $I_0 M_2 I_1 \overline{M_3} I_2 M_1 I_3$\\
\hline
\end{tabular}
\end{table}

%\begin{center}
%\textbf{Paths}                                     \hspace{2.9cm} \textbf{Inversions}
%\end{center}
%$1) \hspace{2mm} e_1(1A \hspace{2mm}2A)$  $I_1M_2 I_2 \overline{M_3}I_3 M_1I_4$ \hspace{10mm}$I_0\overline{ %M_1} I_1M_3I_2 \overline{M_2}I_3$ \\

%\noindent $2) \hspace{2mm} e_1(1A \hspace{2mm}2B)$ $I_1M_2I_2 M_3I_3M_1I_4$ \hspace{10mm}$I_0\overline{M_1}  %I_1 \overline{M_3}I_2 \overline{M_2}I_3$ \\

%\noindent $3) \hspace{2mm} e_1(1B \hspace{2mm}2A)$ $I_1M_2 I_2 \overline{M_3}I_3\overline{M_1}$ \hspace{14mm}
%$M_1  I_1 M_3 I_2\overline{ M_2} I_3$ \\

%\noindent $4) \hspace{2mm} e_1(1B \hspace{2mm}2B)$ $I_1 M_2 I_2 M_3 I_3 \overline{M_1}$ \hspace{14mm} $ M_1 I_0 %\overline{M_3} I_1\overline{ M_2}I_2$\\

%\noindent $5) \hspace{2mm} e_2(1A \hspace{2mm}2A)$ $M_1I_1M_3 I_2\overline{M_2} I_3$ \hspace{14mm} $I_0 M_2 %I_1 \overline{M_3} I_2 \overline{M_1}$\\

%\noindent $6) \hspace{2mm} e_2(1A \hspace{2mm}2B)$ $M_1I_1\overline{ M_3} I_2\overline{M_2}I_3$ %\hspace{14mm}$I_0 M_2 I_1 M_3 I_2\overline{M_1}$\\

%\noindent $7) \hspace{2mm} e_2(1B \hspace{2mm}2A)$ $I_1\overline{M_1}I_2\overline{M_3}I_3 \overline{M_2} I_4$ %\hspace{10mm} $I_0M_2 I_1 M_3 I_2 M_1 I_3$\\

%\noindent $8) \hspace{2mm}e_2(1B \hspace{2mm}2B)$ $I_1\overline{M_1}I_2 M_3 I_3 \overline{M_2}I_4$ %\hspace{10mm} $I_0 M_2 I_1 \overline{M_3} I_2 M_1 I_3$\\
The first column in the table above are the possible HPPs ( $\Gamma$) of assembly graph 1221. The second column represent the micronuclear sequences obtained from the HPPs. The third column is $\Gamma^R$. Notice, the sequences of $\Gamma^R$ are complementary to the DNA sequences of $\Gamma$.  Which leaves eight distinct micronuclear sequences of this fixed orientation. For example:\\
The $\Gamma^{R}$  sequence of line 8 is the sequence in  line 1\\
The $\Gamma^{R}$ sequence of line 7 is the sequence in  line 2\\
The $\Gamma^{R}$ sequence of line 5 is the sequence in  line 3\\
The $\Gamma^{R}$ sequence of line 6 is the sequence in  line 4\\
\indent But, we have to find all possible micronuclear sequences for assembly graph 1221, which also considers $\Gamma^{-R}$ and $\Gamma^-$. The other eight distinct micronuclear sequences come from $\Gamma^{-R}$ and $\Gamma^-$ as follows,

\begin{table}[H]
\begin{tabular}{|c|c|c|}
\hline
HPP Orientation & Micronuclear Sequence\\
\hline
\hline
\multirow{4}{*}{$\Gamma^-$} &  $I_0\overline{M_2} I_1M_1I_2 \overline{M_3}I_3$\\
&  $I_0 \overline{M_2}I_1\overline{M_1}I_2\overline{M_3}I_3$\\
&$I_0 \overline{M_2}I_1M_1I_2M_3I_3$\\
&$I_0 \overline{M_2}I_1\overline{M_1}I_2M_3I_3$\\
\hline
\multirow{4}{*}{$\Gamma^{-R}$} &  $I_0M_3 I_1\overline{M_1}I_2M_2I_3$\\
&  $I_0M_3 I_1M_1I_2M_2I_3$\\
&  $I_0\overline{M_3} I_1\overline{M_1}I_2M_2I_3$\\
&  $I_0\overline{M_3} I_1M_1I_2M_2I_3$\\
\hline
\end{tabular}
\end{table}

with a total of sixteen distinct micronuclear sequences of assembly graph 1221.
\newpage

\section{Smoothing of Assembly Graphs}
\indent Smoothings of assembly graphs represent alignment of pointers (MDS pointers) in the micronuclear sequences and the rearrangement. Each vertex, there are  edges incident to the vertex labeled with MDSs and IESs that cross each vertex. When smoothed at the vertex (removing the vertex), the pointers align (See Figure 10). During this process, MDSs align and the in the micronuclear sequence, recombine MDSs and IESs. There are two types of smoothing, \textit{N} (non-parallel) smoothing and \textit{P} (parallel) smoothing (See Figure 11). If the pointer sequences are aligned in parallel, then the smoothing at the vertex follows the predetermined direction of the assembly graph (parallel smoothing). If the pointer alignment is anti parallel, then the smoothing of the vertex changes the path of orientation of the original graph \cite{1}. Smoothing of each vertex in an assembly graph, is determined by the HPPs of the assembly graph. The type of smoothing is predetermined by the given HPP, in such a way the HPP stays connected. The writing of the "new" micronuclear sequences after the smoothing, is similar to writing of the micronuclear sequence before the smoothing. Use similar notation from Sec.3 to write "new" micronuclear sequences. The difference from the original transversal "old" micronuclear sequence to the new assembly graph "new" micronuclear sequence (not necessary a transversal) is recombined MDSs, since that's what smoothing does to the sequences. Notice, the appearance during P smoothings, the assembly graph is broken up into components, that is why there are bracket notations. The brackets show a separation in MDSs within the assembly graph.\\
\begin{figure}[H]
\includegraphics[width=4in]{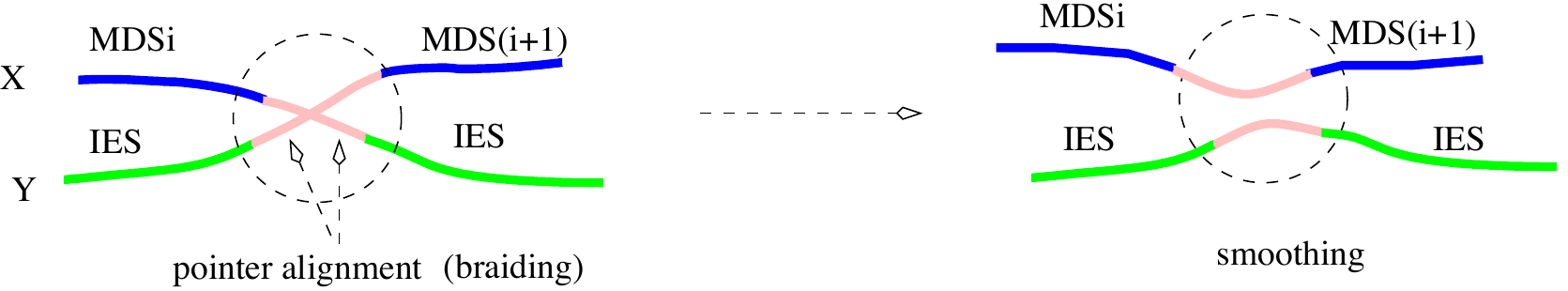}\\
\centering
Figure 10
\end{figure}
\begin{figure}[H]
\includegraphics[width=4in]{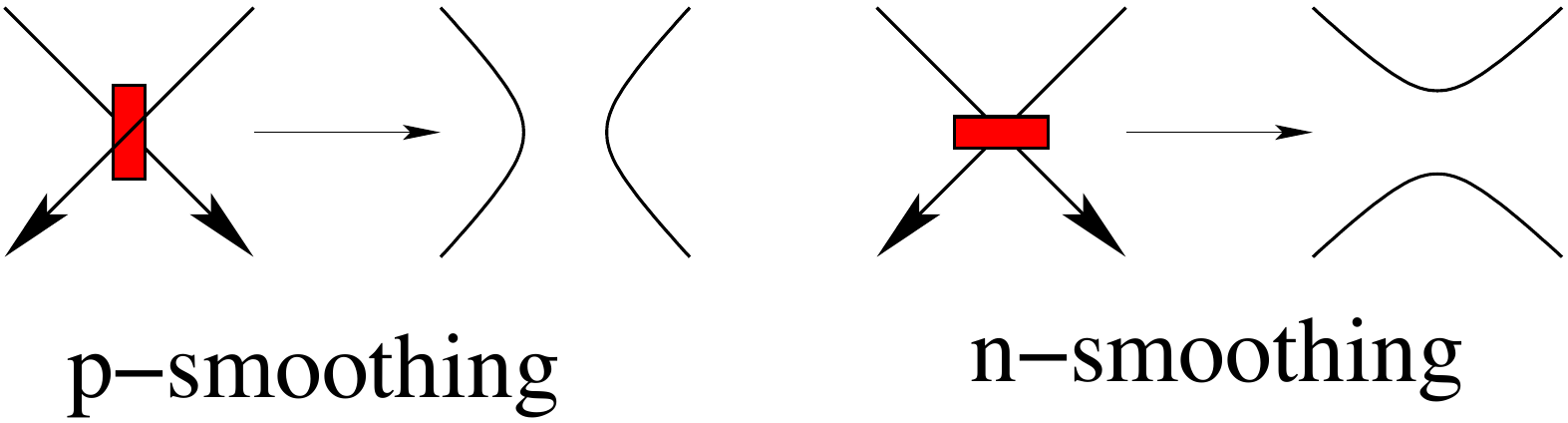}\\
\centering
Figure 11
\end{figure}
\newpage
\subsection{Smoothing of Assembly Graph 1212} 
\indent The following figures represent smoothings of HPPs on assembly graph 1212. Where corresponding micronuclear sequences of the HPPs is on left side of arrow. The smoothings now correspond to a new micronuclear sequence, which is represented on the right side of the arrow. As stated previously, the HPPs predetermines the type of smoothing of the HPPs on the assembly graph, which is one reason why the sequences are grouped in this manner. The other is with respect to which vertex is removed. Notice, the merging of MDSs after vertex removal.\\
\begin{figure}[h]
\begin{minipage}{.45\textwidth}
\begin{center}
\textbf{N smoothing at vertex $1$}\\
$$
\begin{array}{llll}
 1) &I_0 M_1I_1\overline{M_2}I_2 \overline{M_3}I_3 & \longrightarrow & I_0M_{1,2}I_1 \overline{M_3}I_2\\
 2) &I_0  M_1I_1\overline{M_2}I_2 M_3I_3 & \longrightarrow &I_0 M_{1,2}I_1 M_3I_2\\
 7) & I_0\overline{M_1}I_1\overline{M_3}I_2 M_2 I_3 & \longrightarrow & I_0M_{3,1,2}I_1\\
 8) & I_0\overline{M_1}I_1 M_3 I_2 M_2I_3 & \longrightarrow & I_0\overline{M_3}I_1 M_{1,2} I_2\\ 
11) & I_0 M_2 I_1\overline M_1 I_2 \overline{M_3}I_3 & \longrightarrow & I_0\overline M_{2,1,3} I_1\\
12) & I_0 M_2 I_1\overline M_1 I_2 M_3 & \longrightarrow & I_0\overline M_{2,1}I_1 M_3I_2
\end{array}
$$
%\noindent 1)\\
 \begin{figure}[H]
    \includegraphics[width=2.2in] {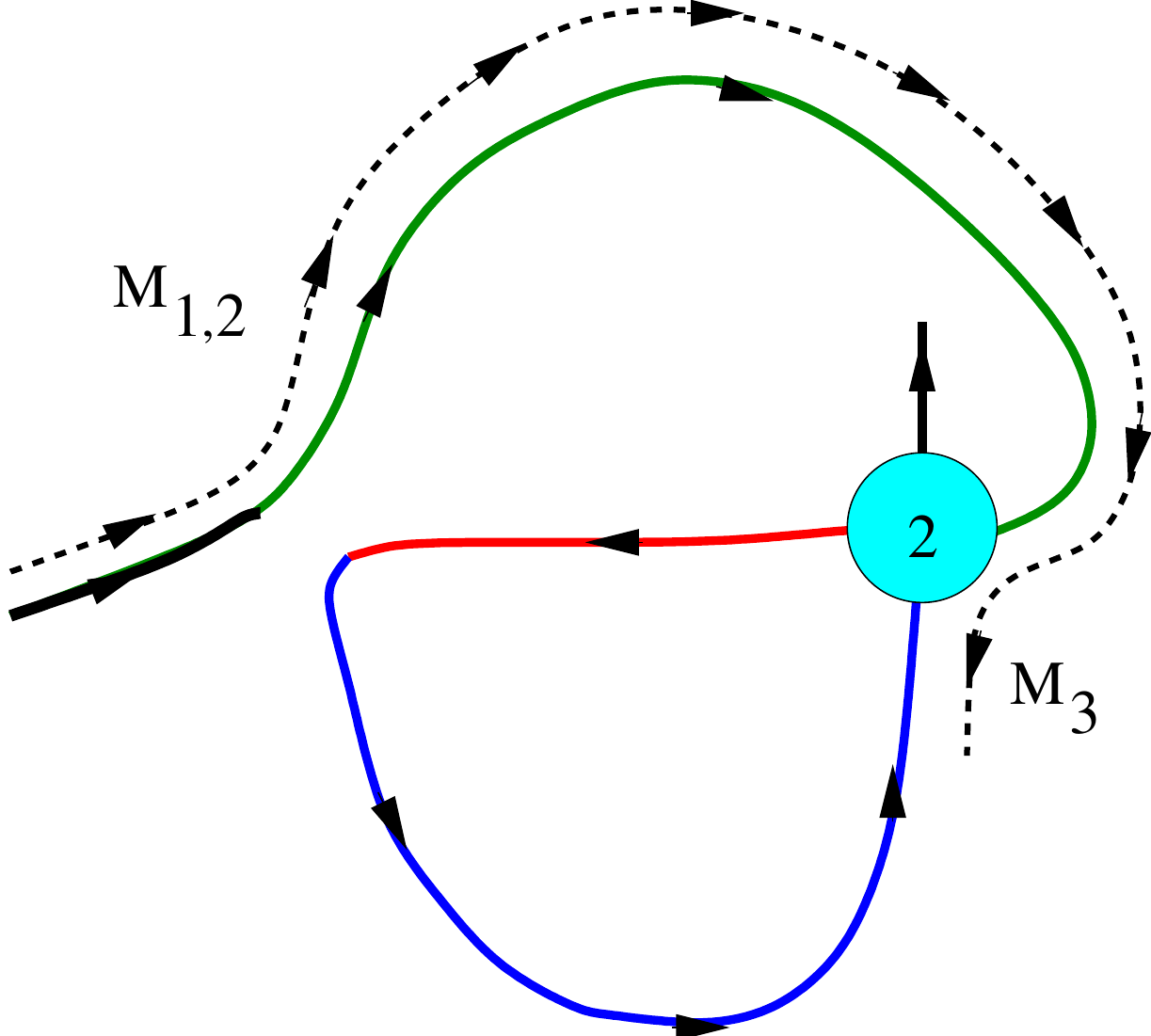}
\caption{The following smoothing corresponds to the first micronuclear sequence of assembly graph 1212.}
\end{figure}
\end{center}
\end{minipage}
\hfill
\begin{minipage}{.45\textwidth}
   \begin{center}
   \noindent \textbf{\noindent P smoothing at vertex $1$\\}
   $$
   \begin{array}{llll}
3)  & I_0\overline{M_1}I_1\overline{M_2}I_2M_3 & \longrightarrow & {I_0M_3,[I_1 \overline{M_{2,1}} I_2]} \\
4)  & I_0\overline{M_1}I_1\overline{M_2}I_2\overline {M_3}I_3 & \longrightarrow & {I_0\overline{M_3} I_1,[\overline{M_{2,1}} I_2]}  \\
5 ) & I_0M_1I_1\overline{M_3}I_2 M_2I_3 & \longrightarrow & I_0M_{1,2}I_1[\overline{M_3}I_2]\\
6)  & I_0M_1I_1 M_3 I_2 M_2I_3 & \longrightarrow &I_0 M_{1,2} I_1[ M_3 I_2]\\
 9)  & I_0 M_2 I_1 M_1 I_2 \overline{M_3}I_3 & \longrightarrow & I_0\overline{M_3} I_1[ M_{1,2}I_2]\\
10)  & I_0 M_2 I_1 {M_1} I_2 {M_3}I_3 & \longrightarrow & I_0 M_3[ I_1 M_{1,2}I_2]
\end{array}
$$
%1)
 \begin{figure}[H]
    \includegraphics[width=2.2in] {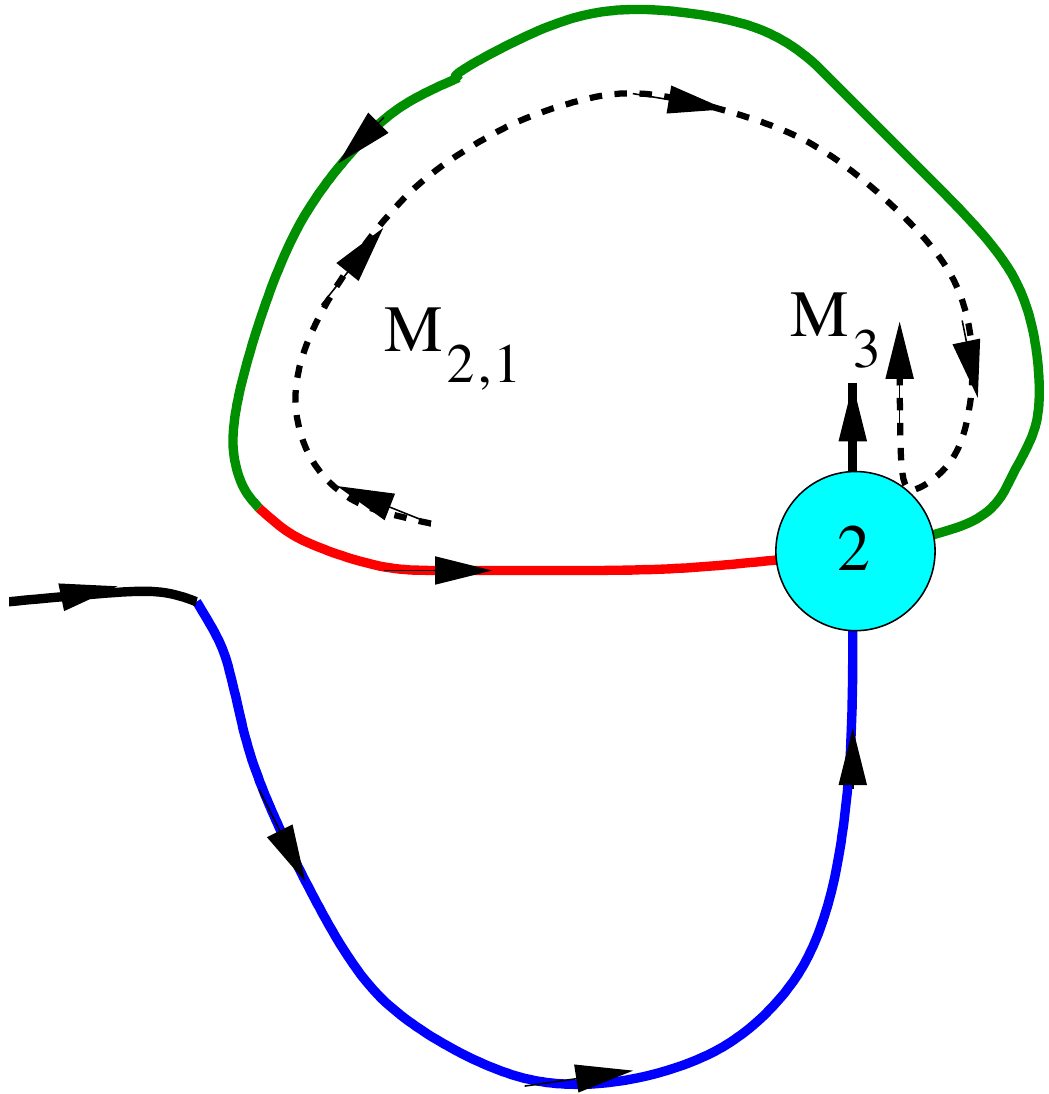}
\caption{The following smoothing corresponds to the third  micronuclear sequence of assembly graph 1212.}
\end{figure}
    \end{center}
\end{minipage}
\end{figure}\\
\indent As seen in Figure 3, represents the assembly graph for first micronuclear sequence of assembly 1212, when vertex one is removed. We can see in the micronuclear sequence before the vertex removal, had separate MDSs. After, the vertex removal, there were merging of MDSs, $M_1$ and $M_2$ to $M_{1,2}$. To obtain the micronuclear sequence, read the transversal path with respect to the smoothed HPP, considering the fixed orientations of both.
\newpage
\begin{figure}[h]
\begin{minipage}{.45\textwidth}
\begin{center}
\noindent \textbf{\noindent P smoothing at vertex $2$\\}
$$
\begin{array}{llll}
1) & I_0M_1I_1\overline{M_2}I_2 \overline{M_3}I_3 & \longrightarrow & I_0M_1I_1[\overline{M_{2,3}}I_2]\\
4)  & I_0\overline{M_1}I_1\overline{M_2}I_2\overline {M_3}I_3 & \longrightarrow & I_0\overline{M_1} I_1[\overline{M_{2,3}}I_2]\\
6) & I_0M_1I_1 M_3 I_2M_2I_3 & \longrightarrow & I_0M_1I_1[M_{3,2}I_2]\\
 8)  & I_0\overline{M_1}I_1 M_3 I_2 M_2I_3 & \longrightarrow & I_0\overline{M_1} I_1 [M_{3,2}I_2]\\
10)  & I_0 M_2 I_1 {M_1} I_2{M_3} I_3& \longrightarrow & I_0M_{2,3}[I_1M_1I_2] \\
12)  & I_0 M_2 I_1\overline M_1 I_2 M_3I_3 & \longrightarrow & I_0M_{2,3}[I_1\overline{M_1}I_2]\\
\end{array}
$$
%3) 
\begin{figure}[H]
    \includegraphics[width=2.2in] {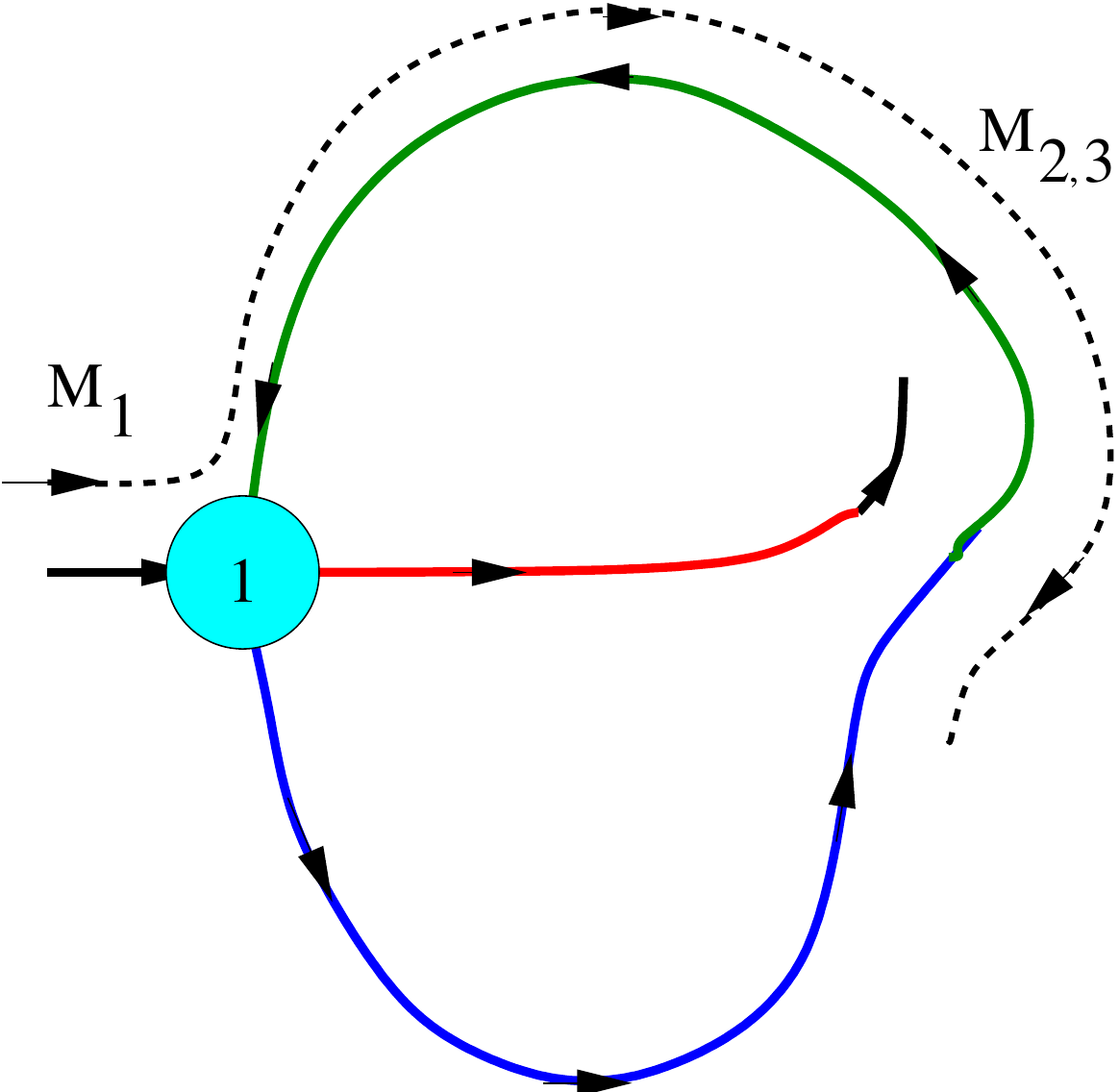}
\caption{The following smoothing corresponds to the first micronuclear sequence of assembly graph 1212.}
\end{figure}
\end{center}
\end{minipage}
\hfill
\begin{minipage}{.45\textwidth}
   \begin{center}
   \noindent \textbf{\noindent N smoothing at vertex $2$\\}
   $$
   \begin{array}{llll}
2)  & I_0M_1I_1\overline{M_2}I_2 M_3I_3 & \longrightarrow  & I_0M_1I_1M_{2,3}I_2\\
3)  & I_0\overline{M_1}I_1\overline{M_2}I_2M_3I_3 & \longrightarrow & I_0\overline{M_1}I_1M_{2,3}I_2\\
5)  &I_0 M_1I_1\overline{M_3}I_2 M_2I_3 & \longrightarrow & I_0M_1I_1\overline{M_{3,2}}I_2\\
7)  & I_0\overline{M_1}I_1\overline{M_3}I_2 M_2 I_3 & \longrightarrow & I_0 \overline{M_{1,3,2}} I_1\\
9)  & I_0 M_2 I_1 M_1 I_2 \overline{M_3}I_3 & \longrightarrow & I_0 M_{2,3}I_1 \overline{M_1}I_2\\
11)  & I_0 M_2 I_1\overline M_1 I_2 \overline{M_3}I_3 & \longrightarrow & I_0 M_{2,3,1}I_1
\end{array}
$$
%2) 
\begin{figure}[H]
    \includegraphics[width=2.2in] {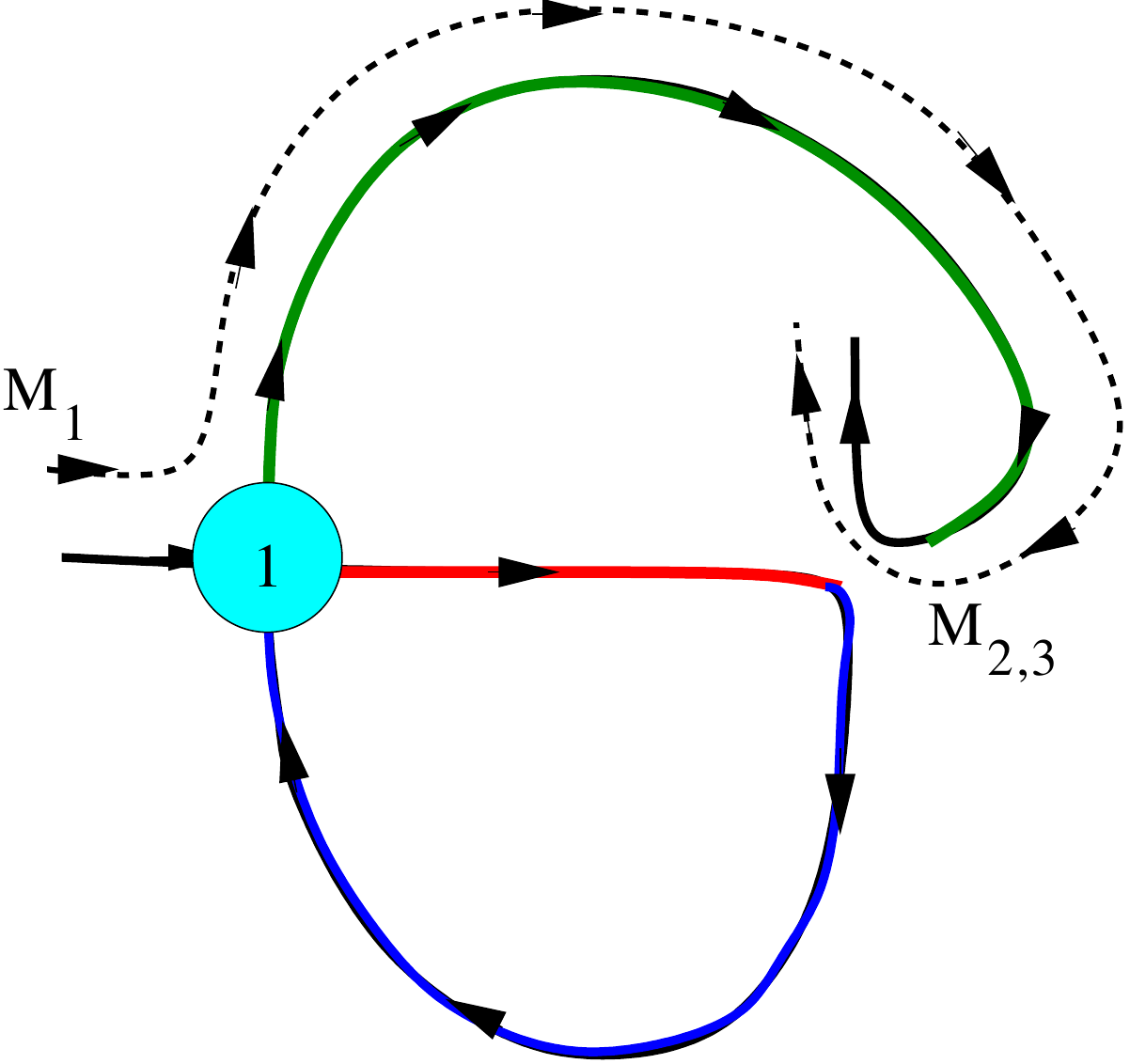}
\caption{The following smoothing corresponds to the second micronuclear sequence of assembly graph 1212.}
\end{figure}
   \end{center}
\end{minipage}
\end{figure}
\indent As seen in Figure 5, represents the assembly graph for first micronuclear sequence of assembly 1212, when vertex two is removed. We can see in the micronuclear sequence before the vertex removal, had separate MDSs. After, the vertex removal, there were merging of MDSs, $\overline{M_2}$ and $\overline{M_3}$ to $\overline{M_{2,3}}$. Notice, there are 2 components of the transversal. To obtain the micronuclear sequence, read the transversal path with respect to the smoothed HPP, considering the fixed orientations of both. The brackets are used to indicate the beginning of the second component of the transversal, with respect to the "new" micronuclear sequence.
\subsection{Smoothing of Assembly Graph 1221}
\indent The following figures represent smoothings of HPPs on assembly graph 1221. Where corresponding micronuclear sequences of the HPPs is on left side of arrow. The smoothings now correspond to a new micronuclear sequence, which is represented on the right side of the arrow. As stated previously, the HPPs predetermines the type of smoothing of the HPPs on the assembly graph, which is one reason why the sequences are grouped in this manner. The other is with respect to which vertex is removed. Notice, the merging of MDSs after vertex removal. Obtaining the "new" micronuclear sequences is similar to smoothings of assembly graph 1212.  \\

%%%%%%   minipage stuff begins
\begin{figure}[h!]
\begin{minipage}{.45\textwidth}
\begin{center}
\textbf{P smoothing at vertex $1$\\}
$$
\begin{array}{llll}
1) &I_0M_2 I_1 \overline{M_3}I_2 M_1I_3 & \longrightarrow  &  I_0M_{1,2}I_1\overline{M_3}I_2\\
2) &I_0M_2I_1 M_3I_2M_1I_3 & \longrightarrow & I_0M_{1,2}I_1M_3I_2\\
7) &I_0\overline{M_1}I_1\overline{M_3}I_2 \overline{M_2} I_3 & \longrightarrow &  I_0\overline{M_{2,1,3}}I_1\\
8) &I_0\overline{M_1}I_1 M_3 I_2 \overline{M_2}I_3 & \longrightarrow & I_0\overline{M_{2,1}}I_1M_3I_2\\
\end{array}
$$
 \begin{figure}[H]
    \includegraphics[width=2.5in] {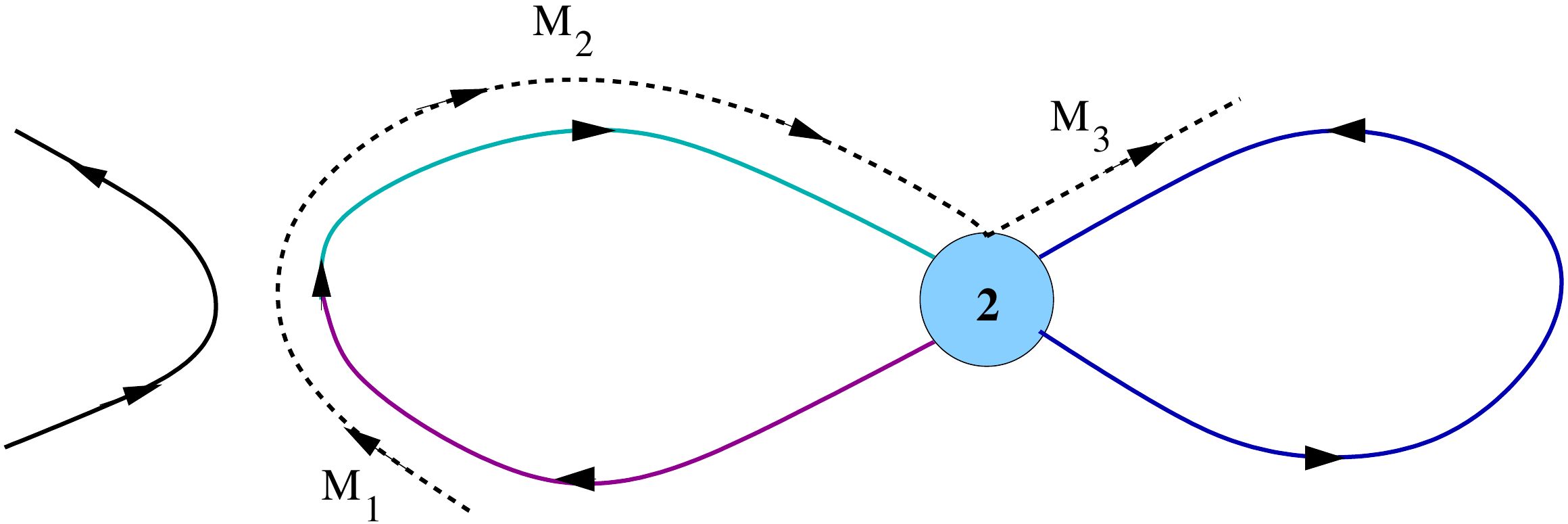}
\caption{{\small The following smoothing corresponds to the first micronuclear sequence of assembly graph 1221.}}
\end{figure}
\end{center}
\end{minipage}
\hfill
\begin{minipage}{.45\textwidth}
   \begin{center}
  \noindent \textbf{N smoothing at vertex $1$\\}  
   $$
   \begin{array}{llll}
 3)  & I_0M_2I_1\overline{ M_3}I_2\overline{M_1}I_3 & \longrightarrow & I_0 M_3 I_1\overline{M_{2,1}}I_2\\
 4)  & I_0 M_2 I_1 M_3 I_2 \overline{M_1}I_0 & \longrightarrow & I_0\overline{M_3} I_1\overline{M_{2,1}}I_2\\
 5) & I_0M_1I_1M_3 I_2\overline{M_2} I_3 & \longrightarrow & I_0M_{1,2}I_1\overline{M_3}I_2\\
 6) & I_0M_1I_1\overline{ M_3} I_2\overline{M_2}I_3 & \longrightarrow & I_0 M_{1,2}I_1M_3I_2
\end{array}
$$
 \begin{figure}[H]
    \includegraphics[width=2in] {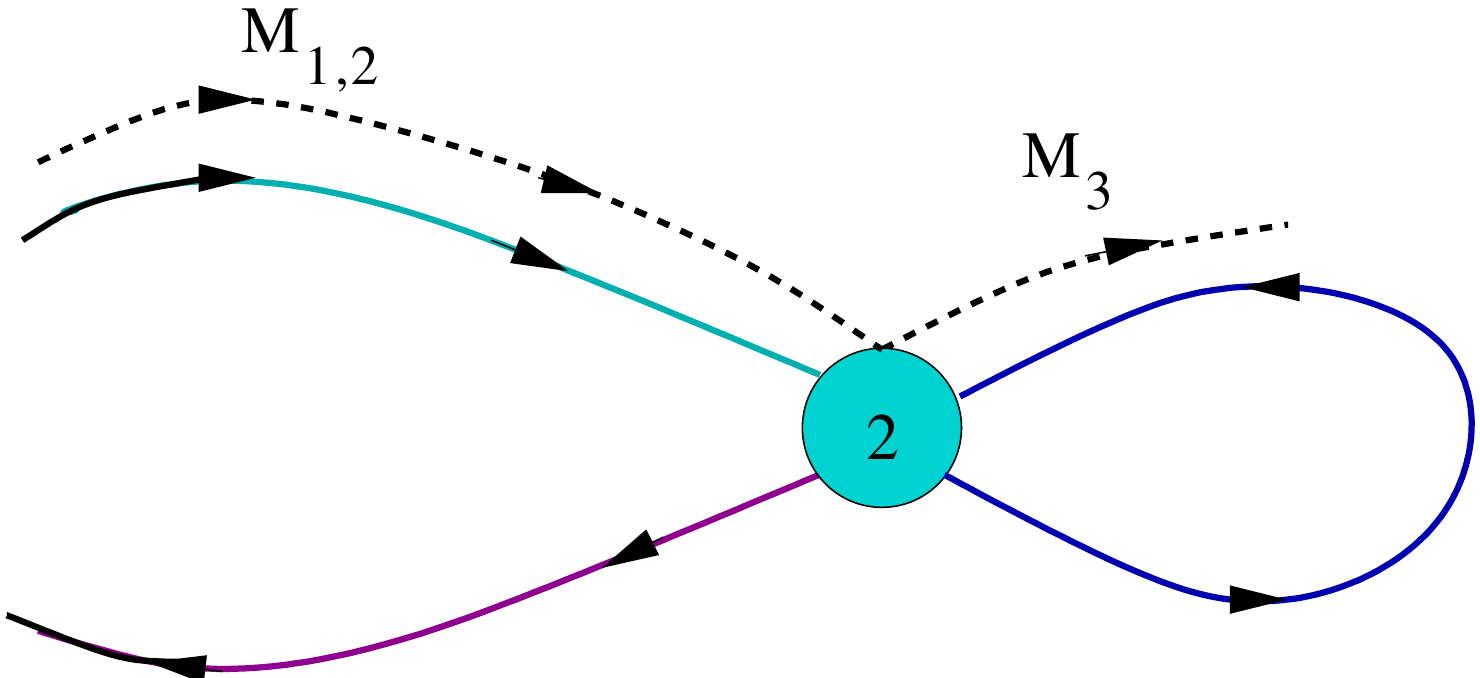}
    \caption{The following smoothing corresponds to the third micronuclear sequence of assembly graph 1221. Removing vertex 1 allowed $\overline{M_1}$ and $M_2$ to align forming $\overline {M_2,1}.$}
\end{figure}
   \end{center}
\end{minipage}
\end{figure}
%%%%%%%%%%%%%%%%%%   end minipage stuff 

%%%%   stuff from here starts in the first minipage part
%
%  \textbf{P smoothing at vertex $1$\\}
%  1) $I_0M_2 I_1 \overline{M_3}I_2 M_1I_3$ $\longrightarrow$ $ I_0M_{1,2}I_1\overline{M_3}I_2$\
%   \
%   2) $I_0M_2I_1 M_3I_2M_1I_3$ $\longrightarrow$ $I_0M_{1,2}I_1M_3I_2$\\
%   7) $I_0\overline{M_1}I_1\overline{M_3}I_2 \overline{M_2} I_3$ $\longrightarrow$ 
%   $I_0\overline{M_{2,1,3}}I_1$\\
%   8) $I_0\overline{M_1}I_1 M_3 I_2 \overline{M_2}I_3$ $\longrightarrow$ 
%   $I_0\overline{M_{2,1}}I_1M_3I_2$\\
% 
%  \begin{figure}[H]
 %  1)  \includegraphics[width=2in] {Smooothing_v1_1}
% \caption{The following smoothing corresponds to the first micronuclear sequence of assembly 
%  graph 1221.}
% \end{figure}
%
%%%%%   first minipage part ends here

%%%%%%% stuff from here to  second part of minipage starts here
%
%\newpage
%\noindent \textbf{N smoothing at vertex $1$\\}
%3) $I_0M_2I_1\overline{ M_3}I_2\overline{M_1}$ $\longrightarrow$ $I_0 M_3 I_1\overline{M_{2,1}}
%$\\
%4) $I_0 M_2 I_1 M_3 I_2 \overline{M_1}$ $\longrightarrow$ $I_0\overline{M_3} 
%I_1\overline{M_{2,1}}$\\
%5) $M_1I_0M_3 I_1\overline{M_2} I_2$ $\longrightarrow$ $M_{1,2}I_0\overline{M_3}I_1$\\
%6) $M_1I_0\overline{ M_3} I_1\overline{M_2}I_2$ $\longrightarrow$ $M_{1,2}I_0M_3I_1$\\

 %\begin{figure}[H]
%  3)  \includegraphics[width=2in] {Smoothing_v1_3}
%\end{figure}
%
%%%%%%%  second part of minipage ends here 

\begin{figure}[H]
\begin{minipage}{.45\textwidth}
\begin{center}
\noindent\textbf{P smoothing at vertex $2$\\}
$$
\begin{array}{llll}
2) & I_0M_2I_1 M_3I_2M_1I_3 & \longrightarrow & I_0M_{2,3} I_1 M_1 I_2\\
4)  & I_0 M_2 I_1 M_3 I_2 \overline{M_1}I_3 & \longrightarrow & I_0 M_{2,3}I_1 \overline{M_1}I_2\\
6)  & I_0M_1I_1\overline{ M_3} I_2\overline{M_2}I_3 & \longrightarrow & I_0M_1I_1\overline {M_{3,2}} 
I_2\\
7)  & I_0\overline{M_1}I_1\overline{M_3}I_2 \overline{M_2} I_3 & \longrightarrow & I_0 
\overline{M_1}I_1\overline{M_{3,2}}I_2
\end{array}
$$
\begin{figure}[H]
    \includegraphics[width=2.5in] {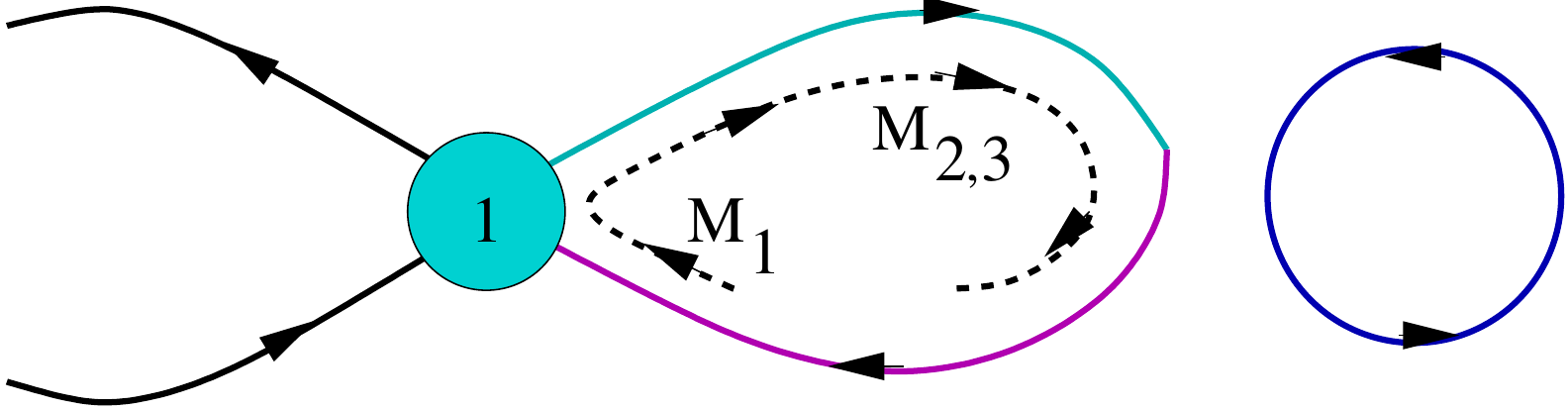}
    \caption {{\small The following smoothing corresponds to the second micronuclear sequence of assembly graph 1221.}}
\end{figure}
\end{center}
\end{minipage}
\hfill
\begin{minipage}{.45\textwidth}
   \begin{center}
   \noindent \textbf{N smoothing at vertex $2$\\}
$$
\begin{array}{llll}
1)  & I_0M_2 I_1 \overline{M_3}I_2 M_1I_3 & \longrightarrow & I_0 M_{2,3}I_1 M_1 I_2\\
3)  & I_0M_2I_1\overline{ M_3}I_2\overline{M_1}I_3 & \longrightarrow & 
I_0 M_{2,3}I_1 \overline{M_1}I_2\\
5)  & I_0M_1I_1M_3 I_2\overline{M_2} I_3 & \longrightarrow & I_0M_1I_1 \overline{M_{3,2}}I_2\\
8)  & I_0\overline{M_1}I_1 M_3 I_2 \overline{M_2}I_3 & \longrightarrow & I_0\overline{M_1}I_1 \overline{M_{3,2}}I_2
\end{array}
$$
 \begin{figure}[H]
    \includegraphics[width=2in] {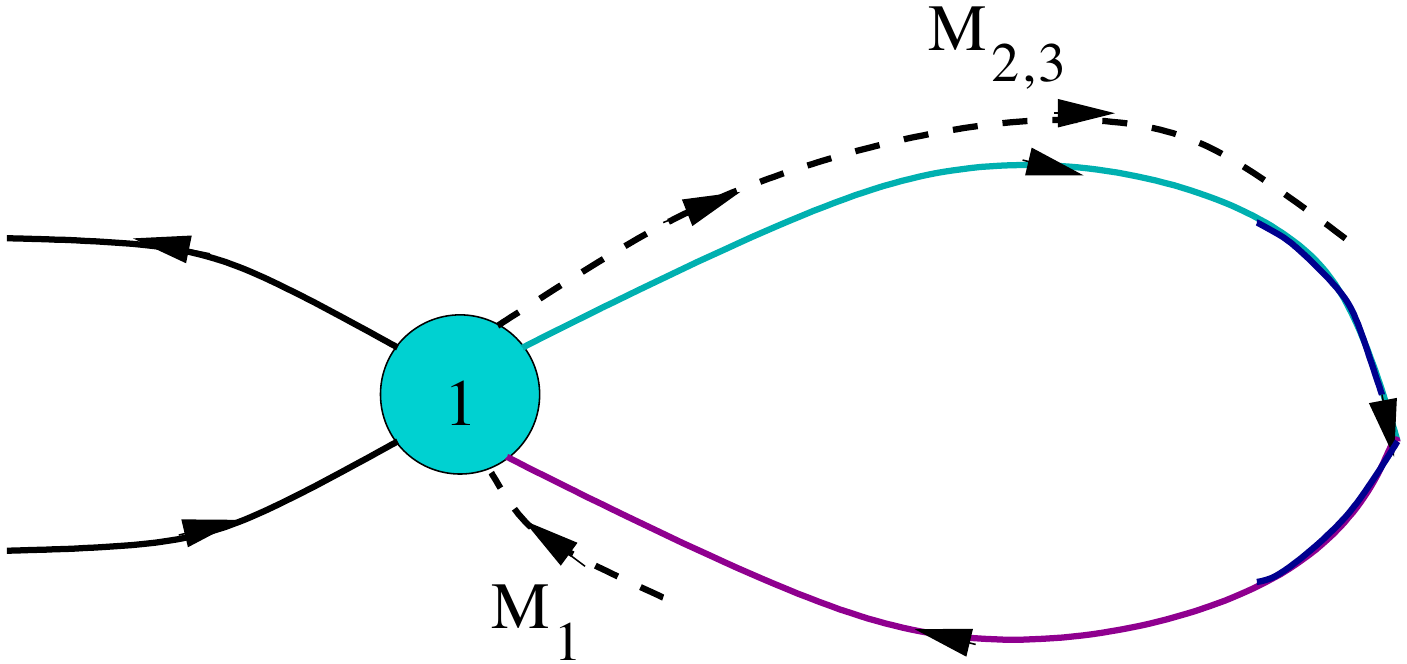}
   \caption{{\small The following smoothing corresponds to the second micronuclear sequence of assembly graph 1221.}}
\end{figure}
   \end{center}
\end{minipage}
\end{figure}

\section{Conclusion}
The forming of micronuclear sequences through HPP can influence the resulting sequences. Such as, the smoothing of vertices, which behave in a manner to preserve the HPP. But, at the same time eliminate some IESs and align MDSs. The previous figures were for two vertices, as the vertices increases, the HPP gets more complicated. The micronuclear sequences are longer, due to the fact more edges are used to form the HPP. The smoothings are more abstract; since the assembly graphs are larger and more smoothings since there are more vertices. There must be four distinct micronuclear sequences for a given path $P(P(1A,2A),...)$. Also, given an orientation of the transversal, each sequence has an inverse. Notice, the second table on page 10 did not include inverses, since they corresponded with the micronuclear sequences of the HPP orientations. As an observation, there are at least eight micronuclear sequence (unless they are complementary or assembly graph 11). For the other orientations $\Gamma^-$ and $\Gamma^{-R}$ there are eight new sequences, potentially.
\section{Acknowledgments}
This work has been supported in part by the NSF grant 
and DMS-0900671. I would like to thank the professors, graduate students and fellow peers who have guided me throughout this summer project. It has truly been an experience long awaited.
 
\end{document}